\RequirePackage{ifpdf}
\ifpdf 
\documentclass[pdftex]{sigma}
\else
\documentclass{sigma}
\fi

\input{xy}
\xyoption{all}

\usepackage{bm}

\def\g{\gamma}
\def\G{\Gamma}

\def\k{\kappa}

\def\L{\Lambda}
\def\O{\Omega}
\def\f{\phi}

\def\o{\omega}
\def\O{\Omega}

\def\s{\sigma}
\def\S{\Sigma}
\def\t{\tau}

\def\CV{{\mathcal{V}}}

\def\P #1{\partial_{#1}}
\def\ch #1{\text{\rm Char}\ #1}

\begin{document}

\allowdisplaybreaks

\renewcommand{\thefootnote}{$\star$}

\renewcommand{\PaperNumber}{098}

\FirstPageHeading

\ShortArticleName{Contact Geometry of Curves}

\ArticleName{Contact Geometry of Curves\footnote{This paper is a
contribution to the Special Issue ``\'Elie Cartan and Dif\/ferential Geometry''. The
full collection is available at
\textit{}\href{http://www.emis.de/journals/SIGMA/Cartan.html}{http://www.emis.de/journals/SIGMA/Cartan.html}}}

\Author{Peter J. VASSILIOU}

\AuthorNameForHeading{P.J. Vassiliou}

\Address{Faculty of Information Sciences and Engineering,
 University of Canberra, 2601 Australia}
\Email{\href{mailto:peter.vassiliou@canberra.edu.au}{peter.vassiliou@canberra.edu.au}}

\ArticleDates{Received May 07, 2009, in f\/inal form October 16, 2009;  Published online October 19, 2009}

\Abstract{Cartan's method of moving frames is brief\/ly recalled in the context of immersed curves in the homogeneous space of a Lie group $G$. The contact geometry of curves in low dimensional equi-af\/f\/ine geometry is then made explicit.  This delivers the complete set of invariant data which solves the $G$-equivalence problem via a straightforward procedure, and which is, in some sense a supplement to the equivariant method of Fels and Olver.
Next, the contact geometry of curves in general Riemannian manifolds $(M,g)$ is described. For the special case in which the isometries of $(M,g)$ act transitively, it is shown that the contact geometry provides an explicit algorithmic construction of the dif\/ferential invariants for curves in~$M$. The inputs required for the construction consist only of the metric~$g$ and a parametrisation of structure group $SO(n)$; the group action is not required and no integration is involved. To illustrate the algorithm we explicitly construct complete sets of dif\/ferential invariants for curves in the Poincar\'e half-space $H^3$ and in a family of constant curvature 3-metrics.  It is conjectured that similar results are possible in other Cartan geometries.}


\Classification{53A35; 53A55; 58A15; 58A20; 58A30}

\renewcommand{\thefootnote}{\arabic{footnote}}
\setcounter{footnote}{0}

\vspace{-2mm}

\section{Introduction}

The classical topic of immersed submanifolds in homogeneous spaces via {\it rep\`ere mobile} or {\it mo\-ving frames} is discussed here in the simplest case, that of curves. Several authors have written on the method of rep\`ere mobile, over the years since Cartan's works, such as~\cite{Cartan37}; these include S.S.~Chern~\cite{Chern85}, J.~Favard~\cite{Favard57}, P.A.~Grif\/f\/iths~\cite{Griffiths74},  G.R.~Jensen~\cite{Jensen77},  M.L.~Green~\cite{Green78},  R.~Sulanke~\cite{Sulanke79}, R.~Sharpe~\cite{Sharpe} and M.E.~Fels \& P.J.~Olver~\cite{FelsOlver98,FelsOlver99}.  Some of these authors have the goal of placing Cartan's method on a f\/irm theoretical foundation as well as extending its range of application beyond the classical realm.  More recently,  a reformulation of the method of moving frames, due to Fels and Olver~\cite{FelsOlver98,FelsOlver99} has lead to renewed activity and a great many new applications and perspectives, have arisen (see~\cite{Olver01} and references therein). Whereas Cartan emphasised the construction of canonical Pfaf\/f\/ian systems whose integral manifolds are the Frenet frames along the submanifold, a much more direct approach is favoured in the Fels--Olver formulation and this has a number of signif\/icant advantages. However, in this paper, we shall reconsider the role of Pfaf\/f\/ian systems in the method of moving frames in the light of recent results in the geometry of jet spaces with the principal goal of making the contact geometry of curves more explicit and giving some indication about its possible applications.  Another goal is to provide additional insight into the relationship between Cartan's method of moving frames and the equivariant method of Fels and Olver\footnote{A point we make herein is that the geometry of jet spaces provides a useful mediation between the two approaches.}.

The considerations in this paper were inspired by a paper of Shadwick and Sluis \cite{ShadSluis}, in which the authors observed that many of the Pfaf\/f\/ian systems derived by Cartan admit a~Cartan prolongation which is locally dif\/feomorphic to the contact distribution on jet space $J^k(\mathbb{R},\mathbb{R}^q)$, for some~$k$ and~$q$, thereby explicitly adding contact geometry to Cartan's method of moving frames.  Another way to view the aims of this paper is the further development of the ideas in~\cite{ShadSluis} in relation to moving frames for curves by exploring the application of a recent generalisation~\mbox{\cite{Vassiliou1,Vassiliou2}} of the Goursat normal form from the theory of exterior dif\/ferential systems~\cite{Stormark2000,BC3G} allowing for the explicit determination of dif\/ferential invariants and other invariant data in cases which have not been previously explored in detail. Of particular interest are curves in general Cartan geo\-metries and in this paper we have focused on the {\it Riemannian} case and conjecture that similar results hold for other Cartan geometries.

\looseness=1
We show that given any $n$-dimensional Riemannian manifold $(M^n,\,g)$ then the Pfaf\/f\/ian system whose integral manifolds determine the Frenet frames along curves in $M$ has a Cartan prolongation which can be identif\/ied with the contact system on jet space $J^n(\mathbb{R},\mathbb{R}^{n-1})$. The explicit construction of the identif\/ication requires only dif\/ferentiation. In case the isometries of $(M^n,g)$ act transitively then the construction of the dif\/ferential invariants that settles the equivalence problem for curves up to an isometry dif\/fers from the approach of Fels--Olver in that explicit a priori knowledge of the isometries or even the inf\/initesimal isometries is not required; as in the Fels--Olver method no integration is called for. The inputs for algorithm {\it Riemannian curves} consist only of the metric $g$ and a realisation of the Lie group $SO(n)$.

Moreover, a contention of this paper is that the contact geometry of submanifolds to be described below {\it should be} a fundamental fact and lead to useful points of view that complement and enhance the geometric analysis of submanifolds by existing methods such as Cartan's method of moving frames and the equivariant moving frames method of Fels and Olver.

\looseness=1
The content of this paper is as follows. After brief\/ly recalling the method of moving frames, as practiced by Cartan, we study one of the simplest non-trivial examples: curves in 2-dimensional equi-af\/f\/ine geometry. It is then shown how the (classical) Goursat normal form applies to give the unique dif\/ferential invariant and moving frame, explicitly. This familiar, illustrative example encapsulates the ideas proposed in this paper and is simple enough so that all details can be given. An account of the generalised Goursat normal form \cite{Vassiliou1,Vassiliou2} is then given in the special case of {\it total} prolongations (uniform Goursat bundles) in preparation for the study of immersed curves in higher dimensional Cartan geo\-metries. Section~\ref{section4} illustrates the princip\-les developed in the previous section by applying it to study curves in 3d-equi-af\/f\/ine geo\-metry, computing the complete set of dif\/ferential invariants via the generalised Goursat normal form.

Section~\ref{section5} is devoted to the contact geometry of curves in any Riemannian manifold and contains the main application of the paper. The general method is used to explicitly derive the dif\/ferential invariants for curves in the Poincar\'e half-space $H^3$ and for curves in a family of constant curvature 3-metrics. These invariants do not seem to have appeared in the literature before. The results demonstrate that the contact geometry of submanifolds can of\/fer an alternative path to invariant data for curves besides the Fels--Olver equivariant method and Cartan's method which, in the latter case, relies so much on geometric insight and special tricks for the construction of the Frenet frames\footnote{In Cartan's writings the distinction between the Frenet frames along a submanifold and the exterior dif\/ferential system whose {\it solutions} are the Frenet frames is sometimes blurred. It's best to keep these two notions quite separate since the construction of the latter is algorithmic while that of the former is not.}.

Finally, it should be mentioned that while this paper only explores the case of curves, the contact geometry of higher dimensional submanifolds could be similarly studied, commencing with the well known characterisation of contact systems in any jet space given in~\cite{Bryant79,Yamaguchi82}.

\section{Method of rep\`ere mobile applied to curves}\label{section2}

According to \cite{Jensen77} the general problem treated by Cartan in~\cite{Cartan37} and elsewhere is that of the {\it invariants of submanifolds} in the homogeneous space of a Lie group~$G$ under the action of~$G$. In this section, I will give a very brief description of the method of {\it rep\`ere mobile}, Cartan's principal tool for addressing this type of problem. More complete discussions can be found in the references quoted above such as \cite{Favard57,Griffiths74,Jensen77, Sulanke79}. The exposition given by Cartan in~\cite{Cartan35} is still well worth reading.

Let $G$ be a Lie group and $H\subset G$ a closed subgroup. Then we have the $H$-principal bundle
\[ 
\xymatrix{
H \ar[r]^{\iota} & G \ar[d]^{\pi}\\
&  G/H
}
\]
of left cosets of $H$ in $G$ and we let $M:=G/H$. Map $\pi$ is the natural projection assigning a left-coset $gH$ to each element of $g\in G$. There is a natural left-action of $G$ on $M$: $g\cdot zH=gzH$, for all $g\in G$. Let $x$ be a local coordinate system on $M$. Cartan typically began with a representation of $G$ which could be ``decomposed'' into columns
$e_1,e_2,\ldots,e_r$ of $H\subset G$ and $\boldsymbol{x}$, a column vector whose components are the coordinates $x$ on $M$.

Cartan def\/ines dif\/ferential 1-forms $\o^i,\ 1\leq i\leq r$ by
\begin{equation}\label{semibasic}
d\boldsymbol{x}=\sum_{i=1}^r\o^i\otimes e_i.
\end{equation}
The 1-forms $\o^i$ are semi-basic for $\pi$. Furthermore, we have 1-forms $\o^j_i$ def\/ined by
\begin{equation}\label{connection}
de_i=\sum_{j=1}^r\o^j_i\otimes e_j.
\end{equation}
The 1-forms $\o^i$, $\o^j_i$, $i,j=1,\ldots,r$ are the components of the Maurer--Cartan form $\o$ on $G$; the integral submanifolds of the Pfaf\/f\/ian system
\[
\o^1=0,\ \  \o^2=0,\ \ \ldots, \ \ \o^r=0
\]
foliates $G$ by the left cosets of $H$.

Suppose $f:T\to M$ is an immersion of a manifold $T$ into $M$. Then a {\it moving frame} is a local map $F:T\to G$ such that $f=\pi\circ F$. That is, the moving frame assigns to each point $t\in T$ a coset $f(t)\in G/H$.
With this general set up, Cartan addresses the following problem for submanifolds of $M$. Let
$f_1:T_1\to M$ and $f_2:T_2\to M$ be submanifolds. Find necessary and suf\/f\/icient conditions, in the form of dif\/ferential invariants, such that there is a~local dif\/feomorphism $\mu:T_1\to T_2$ and element $g\in G$ such that
\begin{equation*}
f_2\circ\mu=g\cdot f_1.
\end{equation*}
The `$\circ$' denotes function composition while `$\cdot$' continues to denote the left-action of $G$ on $M$. A~special case of this is the so-called {\it fixed parametrisation} problem where one takes $T_1=T_2=T$ and $\mu$ is the identity on $T$. This {\it congruence problem} is the one that will be studied in this paper.

In case the submanifolds of $M$ are curves, Cartan begins by choosing a codimension 1 subset of the semibasic 1-forms and def\/ines the Pfaf\/f\/ian system
\begin{equation*}
\O:\ \o^2=0,\ \ \o^3=0,\ \ \ldots  ,\ \ \o^r=0.
\end{equation*}
One studies the solutions of $\O$ since these project via $\pi$ down to curves in $G/H$, which are the objects of interest. One way to do this is via the Cartan--K\"ahler theorem~\cite{BC3G}. Accordingly, one computes the exterior derivatives of the $\o^j$, $j=2,\ldots,r$ and appends these ``integrability conditions'' to $\O$ forming the dif\/ferential ideal~$\bar{\O}$ with independence form $\o^1$. This procedure allows one to prove existence of integral curves for $\O$ and provides information about the number of such integral curves. However, this makes no use of the special origin of the 1-forms in~$\O$, arising as they do from the Maurer--Cartan form $\o$ on $G$. As a result of this one can go much further. From the structure equations of $\o$ and the vanishing of the exterior derivatives $d\o^i$ we deduce additional 1-form equations of the form
\[
\o^j_i-p^j_i\o^1=0,
\]
for some functions $p^j_i$ on $G$, which are appended to $\O$ as integrability conditions, thereby forming the new Pfaf\/f\/ian system
\[
\bar{\O}: \ \o^2=0,\ \ \o^3=0, \ \ \ldots ,\ \ \o^r=0,\ \ \o^j_i-p^j_i\o^1=0.
\]
In essence, the method of moving frames consists of using the fact that $H$ acts on the f\/ibres of $G\to G/H$ on the right inducing a transformation of the Maurer--Cartan form $\o$ on $G$,   \cite[Chapter~7]{Spivak79}. Indeed, the transformation
\begin{equation}\label{frameChange}
(\boldsymbol{x},e_1,\ldots,e_r)\mapsto (\boldsymbol{x},e_1,\ldots,e_r)h,\qquad \forall\;  h\in H
\end{equation}
on $G$ induces the transformation
\begin{equation}\label{connTrans}
\o\mapsto \text{Ad}(h^{-1})\o+h^{-1}dh=\widetilde{\o}
\end{equation}
on the Maurer--Cartan form~$\o$. In turn, this induces a transformation on the functions $p^j_i$. To proceed further we recall the notion of a Cartan prolongation.

\begin{definition}
Let $\mathcal{I}$ be a Pfaf\/f\/ian system on manifold $M$ and $\mathfrak{p}:\widehat{M}\to M$ a f\/ibre bundle. A~Pfaf\/f\/ian system $\widehat{\mathcal{I}}$ on $\widehat{M}$ is said to be a {\it Cartan prolongation} of $(M,\mathcal{I})$ if
\begin{enumerate}\itemsep=0pt
\item[1)] $\mathfrak{p}^*\mathcal{I}\subseteq\widehat{\mathcal{I}}$;
\item[2)] for every integral submanifold ${\sigma}: S\to {M}$ of ${\mathcal{I}}$ there is a unique  integral submanifold $\widehat{\sigma} : S\to \widehat{M}$ of $\widehat{\mathcal{I}}$ that projects to $\sigma$; that is,  $\sigma=\mathfrak{p}\circ\widehat{\sigma}$.
\end{enumerate}
We say that $\widehat{\sigma}$ is the {\it Cartan lift} of $\sigma$.
\end{definition}

If we choose to view $(G\times\mathbb{R}^s,\bar{\O})$, where the factor $\mathbb{R}^s$ carries the ``parameters'' $p^j_i$,  as a~Cartan prolongation of $(G,\O)$ then (\ref{frameChange}) induces a reduction of the trivial bundle $G\times\mathbb{R}^s\to G$ by normalising the coordinates $p^j_i$ on the f\/ibres to simple constants like 0 and $\pm 1$.

Once $\bar{\O}$ has been normalised, the process begins again by taking exterior derivatives of the enlarged, normalised Pfaf\/f\/ian system arising from $\bar{\O}$. Each step selects a subgroup $K\subset H$. If the process terminates at $K=\{\text{identity}\}$ of $G$ then the resulting Pfaf\/f\/ian system arising from~$\bar{\O}$ is canonical. The {\it integral submanifolds} of $\bar{\O}$ are the {\it Frenet frames}, $\mathfrak{F}$. Hereafter we shall label this canonical Pfaf\/f\/ian system by the symbol $\O_{\mathfrak{F}}$.

The main assertion made in this paper is that the canonical Pfaf\/f\/ian system ${\O}_{\mathfrak{F}}$ determining each Frenet frame along an immersed curve admits a Cartan prolongation  $\widehat{\O}_{{\mathfrak{F}}}$ on $E:=G\times \mathbb{R}^\nu$  for some $\nu$, so that $(E,\widehat{\O}_{{\mathfrak{F}}})$ is locally dif\/feomorphic to a jet space $(J^k(\mathbb{R},\mathbb{R}^q),\O^k(\mathbb{R},\mathbb{R}^q))$ where
$\O^k(\mathbb{R},\mathbb{R}^q)\subset T^*J^k(\mathbb{R},\mathbb{R}^q)$ is the contact sub-bundle. The coordinates on the $\mathbb{R}^\nu$ factor of $E$ carry the dif\/ferential invariants of the problem. Indeed, the integral manifolds of $\widehat{\O}_{{\mathfrak{F}}}$,
say, $\widehat{\G}:I\to G\times\mathbb{R}^\nu$ project down to the {\it Frenet lifts} $\G$ of curves $\g:I\to M$ as in
\begin{gather}
\begin{split}
\xymatrix{
& E \ar[d]^{\widehat\pi}\\
& G \ar[d]^{\pi}\\
I \ar[uur]^{\widehat\Gamma} \ar[ur]_{\Gamma} \ar[r]_{\gamma} & G/H
}
\end{split}\label{fundamental}
\end{gather}
where $I\subseteq\mathbb{R}$ is an interval. Hereafter, one of our goals is to give examples which demonstrate the assertion made above, namely that the Pfaf\/f\/ian system $\widehat{\O}_\mathfrak{F}$ can be identif\/ied with a contact system. This identif\/ication can be constructed explicitly and provides explicit coordinate formulas for all the invariant data: dif\/ferential invariants, Fels--Olver equivariant moving frames and invariant dif\/ferential forms. In Section~\ref{section5} we will prove that this procedure can be applied to curves in any Riemannian manifold and in that case it is {\it algorithmic}\footnote{In this paper a construction or procedure is said to be {\it algorithmic} if it can be performed only by dif\/ferentiation and ``algebraic operations'', which includes constructing the inverse of a local dif\/feomorphism. However, integration is strickly excluded.}. Importantly, one is not required to explicitly know the group action {\it a priori}. Before this we will work out some pedagogical examples. The f\/irst of these is suf\/f\/iciently low dimensional so that all details can be given.

\subsection{Curves in the equi-af\/f\/ine plane}\label{section2.1}

The goal in this subsection is to provide a simple illustration of the method of moving frames as described in the previous subsection. We will construct the Frenet frame $\mathfrak{F}$ for a plane curve up to equi-af\/f\/ine transformations by constructing the canonical Pfaf\/f\/ian system $\O_\mathfrak{F}$ and the appropriate Cartan prolongation $\widehat{\O}_\mathfrak{F}$, as described above.  Here equi-af\/f\/ine transformations means the standard transitive action of the Lie group $G=SL(n,\mathbb{R})\ltimes\mathbb{R}^n$ on $\mathbb{R}^n$. For plane curves we take $n=2$; the action on $\mathbb{R}^2$ with local coordinates $\xi_1$, $\xi_2$ is
\begin{equation*}
\left[\begin{matrix}
\xi_1\cr
\xi_2
\end{matrix}\right]\mapsto A
\left[\begin{matrix}
\xi_1\cr
\xi_2
\end{matrix}\right]+\left[\begin{matrix}
x\cr
y
\end{matrix}\right],
\end{equation*}
where $A\in SL(2,\mathbb{R})$ and $x,y\in\mathbb{R}$. We identify $\mathbb{R}^2$ with $G/SL(2,\mathbb{R})$ where the elements of $G$ are matrices of the form
\begin{equation*}
g=\left[\begin{matrix}
1 &0&0\cr
x &a & b\cr
y &c & d
\end{matrix}\right]
\end{equation*}
and $ad-bc=1$. We call this homogeneous space the {\it affine plane} and denote it by $\mathbb{A}^2$. For local coordinates on $\mathbb{A}^2$ we take $\boldsymbol{x}$, the f\/irst column of $g\in G$ and $e_1$, $e_2$ are the next two columns of~$g$. Equations (\ref{semibasic}), (\ref{connection}) give the semi-basic forms
\begin{equation*}
\o^1=\beta dx-bdy,\qquad \o^2=-cdx+ady,
\end{equation*}
and connection forms
\begin{equation*}
\o^1_1=\beta da-b\,dc,\qquad \o^1_2=\frac{1}{a}\big(\beta b\,da+db-b^2dc\big),\qquad \o^2_1=a\,dc-c\,da,
\end{equation*}
where $\beta=(1+bc)/a$, where we have chosen a chart on $SL(2,\mathbb{R})$ in which $a\neq 0$; note that $\o^1_1+\o^2_2=0$. It is useful to record the structure equations
\begin{gather}
d\o^1=\o^1\wedge\o^1_1+\o^2\wedge\o^1_2,\nonumber\\
d\o^2=\o^1\wedge\o^2_1-\o^2\wedge\o^1_1,\nonumber\\
d\o^1_1=\o^2_1\wedge\o^1_2,\label{strEqs_aff2}\\
d\o^1_2=-2\o^1_1\wedge\o^1_2,\nonumber\\
d\o^2_1=2\o^1_1\wedge\o^2_1.\nonumber
\end{gather}
Successive adapted frames are integral curves of certain Pfaf\/f\/ian systems which will be denoted by~$\O^i$, $i=1,2,\ldots$. The f\/irst adapted frames for curves in $\mathbb{A}^2$ are integral curves of the Pfaf\/f\/ian system $\O^1$, consisting of the single 1-form equation
\begin{equation*}
\O^1:\    \o^2=0
\end{equation*}
with independence form $\o^1$.  From structure equations (\ref{strEqs_aff2}), we obtain
$0=d\o^2\equiv \o^1\wedge\o^2_1\mod \o^2$ and hence to complete $\O^1$ to a dif\/ferential ideal $\bar{\O}^1$ we extend it by appending the 2-form equation $\o^2_1\wedge \o^1=0$. This equation implies that there is a function $p$ on $G$ such that the 2-form equation can be replaced by $\o^2_1-p\,\o^1=0$, a kind of ``f\/irst integral''.  We extend $\O^1$ by this 1-form equation and rename the extended Pfaf\/f\/ian system $\bar{\O}^1$ to  get
\begin{equation*}
\bar{\O}^1: \  \o^2=0,\ \  \o^2_1-p\,\o^1=0.
\end{equation*}
Note that the reconstituted $\bar{\O}^1$ is no longer a dif\/ferential ideal.

As discussed in Section~\ref{section2}, an element $h\in H$ acts on the frame $[\boldsymbol{x},e_1,e_2]$ over each point $\boldsymbol{x}\in G/H$ on the right inducing the transformation~(\ref{connTrans})
on the Maurer--Cartan form on $G$. This, in turn induces a transformation on the function $p$. The subgroup $H_1\subset H$ that leaves $\O^1$-invariant has representation
\[
\left[\begin{matrix}
        1&0&0\cr
        0&a&b\cr
        0&0&1/a
        \end{matrix}\right].
\]
We obtain
\[
\widetilde{\o}^1=a^{-1}\o^1,\qquad \widetilde{\o}^2=a\o^2,\qquad \widetilde{\o}^2_1=a^2\o^2_1.
\]
Hence
\[
0=\o^2_1-p\,\o^1=a^{-2}\widetilde{\o}^2_1-p\,a\,\widetilde{\o}^1=
a^{-2}(\widetilde{\o}^2_1-p\,a^3\widetilde{\o}^1).
\]
The Pfaf\/f\/ian system $\bar{\O}^1$ is transformed to
\[
\widetilde{\o}^2=0,\qquad \widetilde{\o}^2_1-p\,a^3\widetilde{\o}^1=0.
\]
That is, the function $p$ undergoes the transformation $p\mapsto a^3p$. Accordingly, we can choose $a$ so that $ap^3=1$ and transform $\bar{\O}^1$ to\footnote{We have made a tacit genericity assumption that $p\neq 0$. The case $p=0$ must be considered separately, as in~\cite{Jensen77}. To simplify the exposition we shall continue to make such genericity assumptions in this paper.}
\[
\O^2:\    \o^2=0, \ \ \o^2_1-\o_1=0.
\]
The integral submanifolds of $\O^2$, with independence form $\o^1$ are the ``second order frames" for curves in $\mathbb{A}^2$. The subgroup $H_2\subset H_1\subset H$ that preserves the elements of $\O^2$ is
\[
\left[\begin{matrix}
        1&0&0\cr
        0&1&b\cr
        0&0&1
        \end{matrix}\right].
\]
We must now extend $\O^2$ to a dif\/ferential ideal by computing the exterior derivative of $\o^2_1-\o^1$. From the structure equations we obtain $3\o^1_1\wedge\o^1=0$. As before, there is a function $q$ on $G$ such that the 2-form equation can be replaced by
\[
\o^1_1-q\o^1=0
\]
so that (the reconstituted) $\bar{\O}^2$ is given by the 1-form equations
\begin{equation}\label{Aff2_2ndOrder}
\bar{\O}^2:\   \o^2=0,\ \   \o^2_1-\o^1=0,\ \   \o^1_1-q\o^1=0
\end{equation}
and is no longer a dif\/ferential ideal.

By performing a $H_2$ change of frame the 1-forms in (\ref{Aff2_2ndOrder}) become
\[
\widetilde{\o}^2=\o^2,\qquad \widetilde{\o}^2_1-\widetilde{\o}^1=\o^2_1-\o^1,\qquad \widetilde{\o}^1_1=\o^1_1-b\,\o^2_1.
\]
Thus $\O^2$ is invariant under a $H_2$ change of frame while
\[
0=\o^1_1-q\o^1=\widetilde{\o}^1_1+(b-q)\widetilde{\o}^1.
\]
We can chose $b=q$ to obtain the 1-form equation $\widetilde{\o}^1_1=0$, and giving rise to the f\/inal adapted frame (dropping tildes)
\[
\O^3:\    \o^2=0,\ \ \o^2_1-\o_1=0,\ \ \o^1_1=0
\]
which ``reduces the isotropy group to the identity''. Thus, $\O^3$ is the Pfaf\/f\/ian system $\O_\mathfrak{F}$ and its integral curves are the Frenet lifts $\mathfrak{F}$ of curves in $\mathbb{A}^2$. Computing the exterior derivative of $\o^1_1=0$ we obtain the 2-form equation $\o^1_2\wedge\o^1=0$ and hence there is a function $\k$ on $G$ such that $\o^1_2-\k\o^1=0$.  This time there is no freedom left in our choice of frame that enables $\k$ to be transformed away. The function $\k$ here is intrinsic.  Hence, in this case, the Cartan prolongation we seek is the Pfaf\/f\/ian system $\O_\mathfrak{F}$ augmented by the 1-form equation $\o^1_2-\k\o^1=0$
\[
\widehat{\O}_\mathfrak{F}:\   \o^2=0,\ \ \o^2_1-\o^1=0,\ \ \o^1_1=0,\ \ \o^1_2-\k\o^1=0,
\]
on $G\times \mathbb{R}_\k$. The integral curves of $\widehat{\O}_\mathfrak{F}$ with independence form $\o^1$ determine the unique equi-af\/f\/ine invariant for plane curves.

\begin{theorem}\label{2dequiaffineCongruence}
Let $I\subseteq\mathbb{R}$ be an interval and $\g_i : I\to \mathbb{A}^2$ be two immersed curves in the equi-affine plane, each parametrised by equi-affine arc-length. Then there is an element $g\in G$ such that $\g_2=g\cdot \g_1$ if and only if their Cartan lifts $\widehat{\G}_i:I\to G\times\mathbb{R}$ satisfy
\begin{equation}\label{equiaffineInvariant}
\widehat{\G}_1^*\k=\widehat{\G}_2^*\k
\end{equation}
identically on $I$.
\end{theorem}

\begin{proof} Let us f\/irstly recall that $(G\times\mathbb{R},\widehat{\O}_\mathfrak{F})\to(G,\O_\mathfrak{F})$ is a Cartan prolongation and $\widehat{\G}_i$ is a~Cartan lift of $\G_i,\ i=1,2$;  see diagram (\ref{fundamental}). By def\/inition, the $\G_i$ are the Frenet lifts of $\g_i$ and~$\widehat{\G}_i$ are Cartan lifts of $\g_i$. Finally $\widehat{\G}_i$ are integral submanifolds of $\widehat{\O}_\mathfrak{F}$ on $G\times\mathbb{R}$ and consquently for $i=1,2$ we have
\begin{gather}
\widehat{\G}_i^*\o^2=\widehat{\G}_i^*\o^1_1=0,\qquad
 \widehat{\G}_i^*\o^1_2=\big(\widehat{\G}_i^*\k\big)\,\big(\widehat{\G}_i^*\o^1\big),\qquad
\widehat{\G}_i^*\o^2_1=\widehat{\G}_i^*\o^1.\label{equiaffineCurves}
\end{gather}

Since both curves are parametrised by {\it equi-affine} arc-length $s$ we have
$\widehat{\G}_1^*\o^1=\widehat{\G}_2^*\o^1=ds$. From this and from (\ref{equiaffineCurves}) we deduce that
\begin{equation}\label{equiaffineMaurerCartan}
\G_1^*\O_{\text{\rm MC}}=\G_2^*\O_{\text{\rm MC}},
\end{equation}
where $\O_{\text{\rm MC}}$ is the Maurer--Cartan form on $G$. It follows from the standard theorem about maps into a Lie group  \cite[Chapter~10, Theorem~18]{Spivak70} that there is a f\/ixed element $g\in G$ such that $\g_2=g\cdot\g_1$.

Conversely, if $\g_2=g\cdot\g_1$ for some $g\in G$, then the Frenet lifts $\G_i$ of $\g_i$
satisfy (\ref{equiaffineMaurerCartan}) and are~integral submanifolds of  $\O_\mathfrak{F}$. But since $\widehat{\O}_\mathfrak{F}$ is a Cartan prolongation of $\O_\mathfrak{F}$, there are Cartan lifts~$\widehat{\G}_i$ of $\G_i$ which are integral submanifolds of $\widehat{\O}_\mathfrak{F}$. Equation (\ref{equiaffineInvariant}) follows from this and equation~(\ref{equiaffineMaurerCartan}).
\end{proof}

\begin{remark}
This theorem encapsulates the basic idea of this paper and is proposed as a model for the study of curves in any Cartan geometry. The relationship between Theorem \ref{2dequiaffineCongruence} and diagram (\ref{fundamental}) should be clear.  The idea now is that by the Goursat normal form $\big(G\times\mathbb{R},\, {\widehat{\O}_\mathfrak{F}}^\perp\big)$ is locally dif\/feomorphic to the jet bundle
$\big(J^4(\mathbb{R},\mathbb{R}),\, \mathcal{C}^{(4)}_1\big)$, where $\mathcal{C}^{(4)}_1$ is the contact sub-bundle of $TJ^4(\mathbb{R},\mathbb{R})$. That is, there is a local dif\/feomorphism
$\phi : G\times\mathbb{R}\to J^4(\mathbb{R},\mathbb{R})$ such that $\phi_*{\widehat{\O}_\mathfrak{F}}^\perp=\mathcal{C}^{(4)}_1$. In fact, Theorem~\ref{2dequiaffineCongruence} proves that knowing the dif\/feomorphism $\phi$ explicitly constructs the unique invariant $\k$ for plane equi-af\/f\/ine curves explicitly, namely the equi-af\/f\/ine curvature, as well as the equi-af\/f\/ine arc-length. Explaining this is the goal of the next subsection.
\end{remark}

\subsection{Goursat normal form}\label{section2.2}

The Goursat normal form is a local characterisation of the contact distribution on $J^k(\mathbb{R},\mathbb{R})$ for all $k\geq 1$, which we denote $\mathcal{C}^{(k)}_1$. The original theorem is not due to Goursat who was its populariser. It appears the theorem is originally due, in some form, to E.~von Weber but the statement of it I give below essentially arises from a 1914 work of Cartan. A good reference is~\cite{Stormark2000}. This reference describes an interesting, relevant but largely forgotten work of Vessiot~\cite{Vessiot26}. First we establish some notation and def\/initions.

\subsubsection{The derived f\/lag}
Suppose $M$ is a smooth manifold and $\CV\subset TM$ a smooth sub-bundle of its tangent bundle. The structure tensor is the homomorphism of vector bundles $\delta :\L^2\CV\to TM/\CV$
def\/ined by
\[
\delta(X,Y)=[X,Y]\mod\CV,\qquad \text{for}\ \ X,Y\in\G(M,\CV).
\]
If $\delta$ has constant rank, we def\/ine the
{\it first derived bundle} $\CV^{(1)}$ as the inverse image of $\delta(\L^2\CV)$ under the canonical projection $TM\to TM/\CV$.
Informally,
\[
\CV^{(1)}=\CV+[\CV,\CV].
\]
The derived bundles $\CV^{(i)}$ are def\/ined inductively:
\[
\CV^{(i+1)}=\CV^{(i)}+[\CV^{(i)},\CV^{(i)}]
\]
assuming that at each iteration it def\/ines a vector bundle, in which case we shall say that $\CV$ is {\it regular}.
For regular $\CV$, by dimension reasons, there will be a smallest $k$ for which $\CV^{(k+1)}=\CV^{(k)}$. This $k$ is called the {\it derived length} of $\CV$ and the whole sequence of sub-bundles
\[
\CV\subset\CV^{(1)}\subset\CV^{(2)}\subset\cdots\subset\CV^{(k)}
\]
the {\it derived flag} of $\CV$. We shall denote by $\CV^{(\infty)}$ the smallest integrable sub-bundle containing $\CV$.

\subsubsection{Cauchy bundles}
Let us def\/ine
\[
\s: \ \CV\to\text{Hom}(\CV,TM/\CV)\qquad \text{by}\ \ \s(X)(Y)=\delta(X,Y)
\]
Even if $\CV$ is regular, the homomorphism $\s$ need not have constant rank. If it does, let us write $\ch\CV$ for its kernel. The Jacobi identity shows that $\ch\CV$ is always integrable. It is called the {\it Cauchy bundle} or {\it characteristic bundle} of~$\CV$.
If $\CV$ is regular and each $\CV^{(i)}$ has a Cauchy bundle then, we say that $\CV$ is {\it totally regular}. Then by the {\it derived type} of $\CV$ we shall mean the list $\{\CV^{(i)},\ch\CV^{(i)}\}$ of subundles.

\vskip 5 pt

\begin{theorem} [Goursat normal form]\label{GoursatNF}
Let $\CV\subset TM$ be a smooth, totally regular, rank~$2$ sub-bundle over smooth manifold $M$ such that
\begin{enumerate}\itemsep=0pt
\item[\text{\rm a)}] $\CV^{(\infty)}=TM$;
\item[\text{\rm b)}] $\dim\CV^{(i+1)}=\dim\CV^{(i)}+1$,\ while $\CV^{(i)}\neq TM$.
\end{enumerate}

Then there is a generic subset $\hat{M}\subseteq M$ such that in a neigbourhood of each point of $\hat{M}$ there are local coordinates $x,z_0,z_1,z_2,\ldots z_k$ such that $\CV$ has local expression
\[
\Bigg\{\P x+\sum_{j=1}^k z_j\P {z_{j-1}},\ \P {z_k}\Bigg\},
\]
where $k=\dim M-2$. That is, $\CV$ is locally equivalent to $\mathcal{C}^{(k)}_1$ on $\hat{M}$.
\end{theorem}

A proof can be found in  \cite[pp.~157--159]{Stormark2000}. The proof of a much more general result in which the Goursat normal form is a special case is given in~\cite{Vassiliou2}. The signif\/icance of the latter is that an procedure is provided for constructing the local contact coordinates $x,z_0,\ldots,z_k$. This is procedure {\it Contact B} on page~287 of~\cite{Vassiliou2} with $\rho_k=1$ and $\rho_1=\rho_2=\cdots=\rho_{k-1}=0$; the $\rho_i$ are def\/ined in Section~\ref{section3}. In this special case we have the following.

\bigskip

\centerline{\it Procedure Contact for the Goursat normal form}

\begin{enumerate}\itemsep=0pt
\item[] $\mathbf{INPUT:}$ Sub-bundle $\CV\subset TM$ of derived length $k$ which  satisf\/ies the hypotheses of Theorem \ref{GoursatNF}.
\item[a)] Fix any f\/irst integral of $\ch\CV^{(k-1)}$, denoted $x$, and  any section  $\boldsymbol{Z}$ of $\CV$ such that $\boldsymbol{Z}x=1$.
\item[b)] Def\/ine a distribution $\Pi^k$ as follows:
\[
\Pi^{l+1}=[Z,\Pi^l],\qquad \Pi^1=\ch\CV^{(1)},\qquad 1\leq l\leq k-1.
\]
\item[c)] Let $z_0$ be any invariant of $\Pi^k$ such that $dx\wedge dz_0\neq 0$.
\item[d)] Def\/ine functions $z_1,z_2,\ldots,z_k$ by $z_j=\boldsymbol{Z}z_{j-1}$, $j=1,\ldots,k$.
\item[] $\mathbf{OUTPUT:}$ Functions $x,z_0,z_1,\ldots,z_k$ are contact   coordinates for $\CV$.
\end{enumerate}

The proof of correctness of this procedure is given in \cite{Vassiliou2}.

\subsection{Equi-af\/f\/ine invariants \& the Goursat normal form}\label{section2.3}

We use this procedure to construct the various invariant objects for this geometry. In fact we will construct the Frenet frames by constructing all the integral submanifolds of $\widehat{\O}_\mathfrak{F}$. So we set $\CV:=\widehat{\O}_\mathfrak{F}^\perp$:
\[
\CV=\big\{\P {\o^1}+\P {\o^2_1}+\k\P {\o^1_2},\ \P \k\big\}.
\]
Note that we have adopted the usual convention of denoting the frame dual to
\[
\o^1,\ \ \o^2,\ \ \o^1_1,\ \ \o^1_2,\ \  \o^2_1
\]
by
\[
\P {\o^1}, \ \ \P {\o^2},\ \  \P {\o^1_1},\ \  \P {\o^1_2},\ \ \P {\o^2_1}.
\]
In local coordinates we have
\[
\CV=\Big\{X:=a\P x+c\P y+b\P a+\k a\P b+\frac{1+bc}{a}\P c,  \ \P \k\Big\}.
\]
Calculation verif\/ies that the hypotheses of Theorem \ref{GoursatNF} are met and that the derived length of $\CV$ is $k=4$. Then step a) of {\it Contact} requires that we construct at least one invariant of
\[
\ch\CV^{(3)}=\big\{\P b,\ \ \P \k,\ \ a\P a+c\P c\big\}
\]
which has invariants $x$, $y$, $a/c$. Any one of these can be taken as the ``independent variable''. Since $x,y$ are local coordinates on $G/SL(2,\mathbb{R})$, we shall choose $x$ for this purpose. It then follows that we may take $\boldsymbol{Z}$ to be
\[
\boldsymbol{Z}=\frac{1}{a}X=\P x+\frac{c}{a}\P y+\frac{b}{a}\P a+\k\P b+\frac{1+bc}{a^2}\P c.
\]
Step   b)  requires the construction of $\Pi^4$. We get
\begin{gather*}
 \Pi^1=\{\P \k\},\qquad \Pi^2=\{\P \k,\ \P b\},\qquad \Pi^3=\Big\{\P \k,\ \P b,\
\frac{1}{a}\P a+\frac{c}{a^2}\P c\Big\},\\
 \Pi^4=\Big\{\P \k,\ \P b,\ \frac{1}{a}\P a+\frac{c}{a^2}\P c,\ \P c\Big\}.
\end{gather*}
The invariants of $\Pi^4$ are in fact $x$, $y$ and hence by step  c), we set $z_0=y$ and construct $z_j=\boldsymbol{Z}z_{j-1}$, $1\leq j\leq 4$. We get
\begin{gather*}
z_1=\frac{c}{a},\qquad z_2=\frac{1}{a^3},\qquad z_3=-\frac{3b}{a^5},\qquad z_4=3\frac{5b^2-a^2\k}{a^7}
\end{gather*}
obtaining the equivalence $\phi:U\subset G\times\mathbb{R}\to J^4(\mathbb{R},\mathbb{R})$ def\/ined by
\[
\phi\big(x,y,a,b,c,\k\big)=\big(x,y,ca^{-1},a^{-3},-3ba^{-5},3(5b^2-a^2\k)a^{-7}\big)
\]
between $\CV$ and the contact distribution on $J^4(\mathbb{R},\mathbb{R})$. The inverse of $\phi$ is
\begin{gather*}
\phi^{-1}\big(x,z_0,z_1,z_2,z_3,z_4\big)= \big(x,z_0,z_2^{-1/3},-3^{-1}z_3z_2^{-5/3},z_1z_2^{-1/3},9^{-1}(5z_3^2-3z_2z_4)z_2^{-8/3}\big)\\
\phantom{\phi^{-1}\big(x,z_0,z_1,z_2,z_3,z_4\big)}{} = (x,y,a,b,c,\k).
\end{gather*}
Hence if we express the curve in $\mathbb{A}^2$ as a graph $(x,f(x))$ then we deduce from $\phi^{-1}$ that the equi-af\/f\/ine curvature is the well known expression
\[
\kappa:=\frac{5f'''(x)^2-3f''(x)f''''(x)}{9f''(x)^{8/3}}=\frac{1}{2}\big(f''(x)^{-2/3}\big)''.
\]
We also obtain the unique $G$-invariant 1-form, the equi-af\/f\/ine arc length by pulling back $\o^1$ by~$\phi^{-1}$,
\[
(\phi^{-1})^*\o^1=f''(x)^{1/3} dx
\]
and the Frenet frame
\[
\mathfrak{F}=\left(\begin{matrix}
1&0&0\cr
x&f''(x)^{-1/3}& -3^{-1}f'''(x)f''(x)^{-5/3}\cr
f(x)&f'(x)f''(x)^{-1/3}& (f''(x)^2-3^{-1}f'''(x)f'(x))f''(x)^{-5/3}\cr
\end{matrix}\right)
\]
of the curve $(x,f(x))$ by pulling back an arbitrary element $g\in G$ by $\phi^{-1}$. Of course, we can express everything in terms an arbitrary immersion $(x(t),y(t))$ into $\mathbb{A}^2$, rather than as a graph.

\begin{remark}
Note that the procedure we have just described for the invariant data of curves in $\mathbb{A}^2$ is {\it not} algorithmic; we had to solve dif\/ferential equations to obtain the equivalence $\phi$. In practice, however, we f\/ind that when the contact system is that of $J^k(\mathbb{R},\mathbb{R}^q)$ where $q>1$, this integration can often be avoided. We will illustrate this for curves in $\mathbb{A}^3$ and prove, in Section~\ref{section5}, that it holds for curves in any Riemannian manifold of dimension greater than~2.
\end{remark}

\section{On the generalised Goursat normal form}\label{section3}

To carry out the programme of the previous section for curves immersed in manifolds of dimension greater than two we must be able to characterise the contact distributions on jet spaces $J^k(\mathbb{R},\mathbb{R}^q)$, for all $k,q\geq 1$; the case $q=1$ being the Goursat normal form. In principle this generalisation should include partial prolongations of the contact distribution on $J^1(\mathbb{R},\mathbb{R}^q)$ and such a characterisation exists -- the generalised Goursat normal form~\cite{Vassiliou1,Vassiliou2}. However, so far the full scope of this characterisation has not been required. It turns out that only {\it total} prolongations of the f\/irst order jet space are suf\/f\/icient. Accordingly, we will only brief\/ly review those parts of \cite{Vassiliou1,Vassiliou2} that are needed for the results to be described in this paper.

\subsection{The singular variety}\label{section3.1}

For each $x\in M$, let
\[
\mathcal{S}_x=\{v\in\CV_x\backslash 0~|~\s(v)\ \text{has less than generic rank}\}.
\]
Then $\mathcal{S}_x$ is the zero set of homogeneous polynomials and so def\/ines a subvariety of the projectivisation
$\mathbb{P}\CV_x$ of $\CV_x$. We shall denote by Sing$(\CV)$ the f\/ibre bundle over $M$ with f\/ibre over $x\in M$ equal to
$\mathcal{S}_x$ and we refer to it as the {\it singular variety} of $\CV$. For $X\in \CV$ the matrix of the homomorphism $\s(X)$ will be called the {\it polar matrix} of $[X]\in\mathbb{P}\CV$. There is a map $\text{deg}_\CV:\mathbb{P}\CV\to
\mathbb{N}$ well def\/ined by
\[
\text{deg}_\CV([X])=\text{rank}~\s(X)\qquad \text{for}\ \ [X]\in\mathbb{P}\CV.
\]
We shall call $\text{deg}_\CV([X])$ the {\it degree} of $[X]$. The singular variety $\text{Sing}(\CV)$ is a dif\/feomorphism invariant in the sense that if $\CV_1$, $\CV_2$ are sub-bundles over
$M_1$, $M_2$, respectively and there is a~dif\/feomorphism $\phi: M_1\to M_2$ that identif\/ies them, then
$\text{Sing}(\CV_2)$ and $\text{Sing}(\phi_*\CV_1)$ are equiva\-lent as projective subvarieties of $\mathbb{P}\CV_2$. That is, for each $x\in M_1$, there is an element of the projective linear group $PGL({\CV_2}_{|_{\f(x)}},\mathbb{R})$ that identif\/ies
$\text{Sing}(\CV_2)(\f(x))$ and $\text{Sing}(\phi_*\CV_1)(\f(x))$.

We hasten to point out that the computation of the singular variety for any given sub-bundle  $\CV\subset TM$ is algorithmic. One computes the determinantal variety of the polar matrix for generic~$[X]$.

\subsubsection{The singular variety in positive degree}

If $X\in\ch\CV$ then $\text{deg}_\CV([X])=0$. It is convenient to eliminate lines of degree zero and for this reason we pass to the quotient $\widehat{\CV}:=\CV/\ch\CV$. We have structure tensor
$\widehat{\delta}:\L^2\widehat{\CV}\to \widehat{TM}/\widehat{\CV}$, well def\/ined by
\[
\widehat\delta(\widehat{X},\widehat{Y})=\pi([X,Y])\mod\widehat{\CV},
\]
where $\widehat{TM}=TM/\ch\CV$ and
\[
\pi: \ TM\to \widehat{TM}
\]
is the canonical projection. The notion of degree descends to this quotient giving a map
\[
\text{deg}_{\widehat{\CV}}:\ \mathbb{P}\widehat{\CV}\to\mathbb{N}
\]
well def\/ined by
\[
\text{deg}_{\widehat{\CV}}([\widehat{X}])=\text{rank}~\widehat{\s}(\widehat{X})\qquad \text{for}\ \
[\widehat{X}]\in\mathbb{P}\widehat{\CV},
\]
where $\widehat{\s}(\widehat{X})(\widehat{Y})=\widehat{\delta}(\widehat{X},\widehat{Y})$ for $\widehat{Y}\in\widehat{\CV}$. Note that all def\/initions go over {\it mutatis mutandis} when the structure tensor $\delta$  is replaced by
$\widehat{\delta}$. In particular, we have notions of polar matrix and singular variety, as before. However, if the singular variety of $\widehat{\CV}$ is not empty, then each point of
$\mathbb{P}\widehat{\CV}$ has degree one or more.

\subsubsection{The resolvent bundle}  Suppose $\CV\subset TM$ is totally regular of rank $c+q+1$, $q\geq 2$, $c\geq 0$,
$\dim M=c+2q+1$. Suppose further that $\CV$ satisf\/ies
\vskip 2 pt
\begin{itemize}\itemsep=0pt
\item[\rm a)] $\dim\ch\CV=c$, $\CV^{(1)}=TM$;
\item[\rm b)]
$\widehat{\S}_{|_x}:=\text{Sing}(\widehat{\CV})_{|_x}=\mathbb{P}\widehat{\mathcal{B}}_{|_x}\approx\mathbb{R}
\mathbb{P}^{q-1}$, for each $x\in M$ and some rank $q$ sub-bundle $\widehat{\mathcal{B}}\subset\widehat{\CV}$.
Then we call  $(\CV,\mathbb{P}\widehat{\mathcal{B}})$ (or $(\CV,\widehat{\S})$)  a
{\it   Weber structure} of rank $q$ on $M$.
\end{itemize}

Given a  Weber structure
$(\CV,\mathbb{P}\widehat{\mathcal{B}})$, let $\mathcal{R}(\CV)\subset\CV$, denote the largest sub-bundle such that
\begin{equation}\label{2.1}
{\pi}\big( \mathcal{R}(\CV) \big)= \widehat{\mathcal{B}}.
\end{equation}
We call the rank $q+c$ bundle $\mathcal{R}(\CV)$ def\/ined by \eqref{2.1} the {\it resolvent bundle} associated to the  Weber structure $({\CV},\widehat{\S})$. The bundle $\widehat{\mathcal{B}}$ determined by the singular variety of
$\widehat{\CV}$ will be called the {\it singular sub-bundle} of the  Weber structure. A  Weber structure
will be said to be {\it integrable} if its resolvent bundle is integrable.

We will see that the resolvent bundle is the key to the construction of an identif\/ication of a~given dif\/ferential system with a contact system, if such an identif\/ication exists; and hence the name.

An {\it integrable}  Weber structure descends to the quotient of $M$ by the leaves of $\ch\CV$ to be the contact bundle on $J^1(\mathbb{R},\mathbb{R}^q)$. Thus, the resolvent bundle and its concomitant Weber structure is a constructive characterisation of the contact bundle on the 1-jets $J^1(\mathbb{R},\mathbb{R}^q)$. The term `Weber structure' honours Eduard von Weber (1870--1934) who was the f\/irst to publish a~proof of the Goursat normal form. For completeness we record the following properties of the resolvent bundle of a Weber structure.

\begin{proposition}[\cite{Vassiliou1}]\label{Weber_prop}

Let $(\CV,\widehat{\S})$ be a Weber structure on $M$ and
$\widehat{\mathcal{B}}$ its singular sub-bundle. If $q\geq 3$, then the following are equivalent

\begin{itemize}\itemsep=0pt
\item[\rm a)] its resolvent bundle $ \mathcal{R}(\CV) \subset\CV$ is integrable;
\item[\rm b)] each point of $\widehat{\S}=\text{\rm Sing}(\widehat{\CV})$ has degree one;
\item[\rm c)] the structure tensor $\widehat{\delta}$ of $\widehat{\CV}$ vanishes on $\widehat{\mathcal{B}}$:
$\widehat{\delta}(\widehat{\mathcal{B}},\widehat{\mathcal{B}})=0$.
\end{itemize}
\end{proposition}

\begin{proposition}[\cite{Vassiliou1}]  Let $({\CV},\widehat{\S})$ be an integrable Weber structure on $M$.
Then its resolvent bundle $\mathcal{R}(\CV)$ is the unique, maximal, integrable sub-bundle of $\CV$.
\end{proposition}

Checking the integrability of the resolvent bundle is algorithmic.
One computes the singular variety $\text{Sing}(\widehat{\CV})=\mathbb{P}\widehat{\mathcal{B}}$.
In turn, the singular bundle $\widehat{\mathcal{B}}$ algorithmically determines  $\mathcal{R}(\CV)$.

\begin{example}
The canonical model of an integrable Weber structure is the contact distribution on $J^1(\mathbb{R},\mathbb{R}^q)$ for $q>1$, extended by Cauchy characteristics
\[
\mathcal{V}=\Bigg\{\P x+\sum_{i=1}^q p_i\P {u_i},\ \P {p_1},\ldots,\P {p_q},\ \P {z_1},\,\P {z_2},\ldots,
\P {z_c}\Bigg\}.
\]
The quotient $\widehat{\CV}=\mathcal{\CV}/\ch\CV$ has singular sub-bundle
\[
\widehat{\mathcal{B}}=\big\{\left[\P {p_1}\right],\ldots,\left[\P {p_q}\right]\big\}
\]
and the resolvent bundle  of integrable Weber structure $(\CV,\mathbb{P}\widehat{\mathcal{B}})$ is
\[
\mathcal{R}(\CV)=\big\{\P {p_1},\ldots,\P {p_q},\ \P {z_1},\,\P {z_2},\ldots,
\P {z_c}\big\}.
\]
The invariants of the resolvent bundle are spanned by $\{x,\, u_1,\ldots, u_q\}$. So the resolvent bundle provides a geometric characterisation of the ``independent variable''~$x$ and the ``dependent variab\-les''~$u_i$, after which dif\/ferentiation by a canonically def\/ined total derivative operator leads to higher order jet coordinates. See Sections~\ref{section3.2} and~\ref{section4} of this paper for further details. See \cite{Vassiliou1,Vassiliou2} for the general theory with proofs and further examples.
\end{example}

\subsection{The uniform generalised Goursat normal form}\label{section3.2}

We are now able to give a characterisation of the contact distribution on $J^k(\mathbb{R},\mathbb{R}^q)$, $\mathcal{C}^{(k)}_q$ for any $k,q\geq 1$, generalising the Goursat normal form to the uniform case.

\begin{theorem}[generalised Goursat normal form -- uniform case, \cite{Vassiliou1,Vassiliou2}]\label{genGoursatNF_uniform}
Let $\CV\subset TM$ be a smooth, totally regular, sub-bundle of rank $q+1$ and derived length $k$, some $k,q>0$, over smooth manifold $M$ such that{\samepage
\begin{itemize}\itemsep=0pt
\item[{\rm a)}] $\CV^{(\infty)}=TM$;
\item[{\rm b)}] $\dim\CV^{(i+1)}=\dim\CV^{(i)}+q$, while $\CV^{(i)}\neq TM$;
\item[{\rm c)}] $\ch\CV^{(i)}\subset\CV^{(i-1)}$, $1\leq i\leq k-1$;  $\dim\ch\CV^{(j)}=jq$, $0\leq j\leq k-1$;
\item[{\rm d)}] If $q>1$ then $\CV^{(k-1)}$ admits an integrable Weber structure.
\end{itemize}}
Then there is a generic subset $\hat{M}\subseteq M$ such that around each point of $\hat{M}$ the distribution $\CV$ is locally equivalent to $\mathcal{C}^{(k)}_q$.
\end{theorem}

If $q=1$ then c) follows from b) and the Weber structure is not def\/ined in which case Theorem~\ref{genGoursatNF_uniform}
reduces to the Goursat normal form, Theorem~\ref{GoursatNF}. We call any sub-bundle that satisf\/ies the hypotheses of Theorem~\ref{genGoursatNF_uniform} a {\it uniform Goursat bundle} in which case the Theorem asserts that generically every uniform Goursat bundle is locally equivalent to the canonical one. A proof of Theorem \ref{genGoursatNF_uniform} can be found in~\cite{Vassiliou1} as a special case of that paper's Theorem~4.1. However, the latter covers a very much larger class of sub-bundles than uniform Goursat bundles and there is therefore a much simpler proof in this uniform case. However, for the purposes of this paper an important thing is the procedure for constructing contact coordinates in the uniform case which is a special case of procedure {\it Contact A} on page~286 of \cite{Vassiliou2}\footnote{Unfortunately in reference \cite{Vassiliou2} procedure {\it Contact} was called an algorithm. Manifestly, it does not qualify as an algorithm because in a certain step an integration is required.} with
$\rho_1=\rho_2=\cdots=\rho_{k=1}=0$, $\rho_k=q>1$.  Note that the collection of non-negative integers
\[
\sigma=\langle \rho_1,\rho_2,\ldots,\rho_k\rangle
\]
shall be called the {\it signature} of $\mathcal{V}$ and is a complete local invariant of Goursat bundles. In the interests of completeness we mention that the non-negative integers $\rho_i$ are def\/ined by
\begin{gather*}
 \rho_j=\dim\ch\CV^{(j)}-\dim\ch\CV^{(j)}_{j-1},\qquad 1\leq j\leq k-1,\\
 \rho_k=\dim\CV^{(k)}-\dim\CV^{(k-1)},
\end{gather*}
where $\ch\CV^{(j)}_{j-1}=\ch\CV^{(j)}\cap\CV^{(j-1)}$. It is proved in \cite{Vassiliou1,Vassiliou2} that a sub-bundle $\mathcal{V}$ on manifold~$M$ is locally dif\/feomorphic to a partial prolongation of the contact system $\mathcal{C}^{(1)}_q$ on $J^1(\mathbb{R},\mathbb{R}^q)$ with~$\rho_j$ variables of order $j$, if and only if $(M,\, \mathcal{V})$ is a Goursat bundle of signature~$\sigma$. In this paper, we need only consider total prolongations in which the only nonzero element of the signature is $\rho_k$, where $k$ is the derived length of $\mathcal{V}$.

\bigskip

\centerline{\it Procedure Contact for uniform
 generalised Goursat bundles, $q>1$}

\begin{enumerate}\itemsep=0pt
\item[] $\mathbf{INPUT:}$ Uniform Goursat bundle $\CV\subset TM$ of  derived  length $k$ with $q>1$.
\item[a)] Construct the (integrable) resolvent bundle  $\mathcal{R}(\CV^{(k-1)})$ and all its $q+1$ f\/irst integrals.
\item[b)] Fix any one of the f\/irst integrals from step a)  denoted $x$, and any section $\boldsymbol{Z}$ of $\CV$ such that $\boldsymbol{Z}x=1$.
\item[c)] Denote the remaining $q$ f\/irst integrals of $\mathcal{R}(\CV^{(k-1)})$  by $z_0^j$, $j=1,2,\ldots, q$.
\item[d)] Def\/ine functions $z^j_1,z^j_2,\ldots,z^j_k$ by $z^j_m=\boldsymbol{Z}z^j_{m-1}$,  $1\leq m\leq k$, $1\leq j\leq q$.
\item[] $\mathbf{OUTPUT:}$ Functions $x,z^j_0,z^j_1,\ldots,z^j_k$, $1\leq j\leq q$ are contact coordinates for $\CV$.
\end{enumerate}

\begin{remark}
Even though an integration problem is presented for solution in step~a) above, in fact, in every example of curves in a Cartan geometry that I've seen no integration is required because the resolvent bundle turns out to be the vertical bundle for the f\/ibration $\pi\circ\widehat{\pi}:$ \mbox{$E\to G/H$}. So a complete set of invariants of $\mathcal{R}(\CV^{(k-1)})$ can be taken to be the components of any coordinate system on $G/H$. In Section~\ref{section5} we will prove this for curves in any Riemannian manifold of dimension at least~3.
\end{remark}

\section{Space curves up to equi-af\/f\/ine transformations}\label{section4}

As an illustration of the generalised Goursat normal form and its relation to the geometry of curves we consider immersed curves in $\mathbb{R}^3$ up to the standard action of $G:=SL(3,\mathbb{R})\ltimes\mathbb{R}^3$. The goal
is to use rep\`ere mobile and Theorem \ref{genGoursatNF_uniform} to construct all the invariant data for this situation.

We discuss this example principally for illustration since we permit ourselves to begin with the explicit transitive action. However, our goal in this paper is to drop any reliance on an {\it a priori} knowledge of a group action. To that end a contention of this paper is that explicit invariant curve data {\it can} be obtained without integration or explicit knowledge of the group action in a signif\/icant special case, namely Riemannian geometry. This will be established in Section~\ref{section5}.

A straightforward extension of the $n=2$ case covered in Section~\ref{section2.1} to the $n=3$ case leads to the matrix group with elements
\begin{equation*}
g=\left(\begin{matrix}
1&0&0&0\cr
x&a_1&a_2&a_3\cr
y&a_4&a_5&a_6\cr
z&a_7&a_8&a_9
\end{matrix}\right),
\end{equation*}
where $\det g=1$. We parametrise an open subset of the group $G$ by solving $\det g=1$ for $a_9$. From equations (\ref{semibasic}), (\ref{connection}) or otherwise we deduce the left-invariant Maurer--Cartan form
\begin{equation*}
\omega=
\left(\begin{matrix}
0 & 0 & 0&0\cr
\o^1 & \o^1_1 & \o^1_2&\o^1_3\cr
\o^2 & \o^2_1 & \o^2_2&\o^2_3\cr
\o^3 & \o^3_1 & \o^3_2&\o^3_3\cr
\end{matrix}\right),
\end{equation*}
where $\o^1_1+\o^2_2+\o^3_3=0$.
For a coframe on $G$ we take the ordered list
\[
\big[\o^1,\o^2,\o^3,\o^1_1,\o^2_1,\o^3_1,\o^1_2,\o^2_2,\o^3_2,\o^1_3,\o^2_3\big]
\]
whose dual frame we label
$
\mathfrak{V}=\big[v_1,v_2,\ldots,v_{11}\big].
$
The Lie algebra multiplication table for $\mathfrak{V}$ is

\centerline{\begin{tabular}{c|ccccccccccc}
  &$v_1$&$v_2$&$v_3$&$v_4$&$v_5$&$v_6$&$v_7$&$v_8$&$v_9$&$v_{10}$&$v_{11}$\cr
\cline{1-12}
$v_1$&$0$&$0$&$0$&$-v_1$&$-v_2$&$-v_3$&0&0&0&0&0\cr
$v_2$&$0$&$0$&0&0&0&0&$-v_1$&$-v_2$&$-v_3$&0&0\cr
$v_3$&0&0&0&$v_3$&0&0&0&$v_3$&0&$-v_1$&$-v_2$\cr
$v_4$&$v_1$&0&$-v_3$&0&$-v_5$&$-2v_6$&$v_7$&0&$-v_9$&$2v_{10}$&$v_{11}$\cr
$v_5$&$v_2$&0&0&$v_5$&0&0&$v_8-v_4$&$-v_5$&$-v_6$&$v_{11}$&0\cr
$v_6$&$v_3$&0&0&$2v_6$&0&0&$v_9$&$v_6$&0&$-v_4$&$-v_5$\cr
$v_7$&0&$v_1$&0&$-v_7$&$v_4-v_8$&$-v_9$&0&$v_7$&0&0&$v_{10}$\cr
$v_8$&0&$v_2$&$-v_3$&0&$v_5$&$-v_6$&$-v_7$&0&$-2v_9$&$v_{10}$&$2v_{11}$\cr
$v_9$&0&$v_3$&0&$v_9$&$v_6$&0&0&$2v_9$&0&$-v_7$&$-v_8$\cr
$v_{10}$&0&0&$v_1$&$-2v_{10}$&$-v_{11}$&$v_4$&0&$-v_{10}$&$v_7$&0&0\cr
$v_{11}$&0&0&$v_2$&$-v_{11}$&0&$v_5$&$-v_{10}$&$-2v_{11}$&$v_8$&0&0\cr
\end{tabular}}

\medskip

By a procedure similar to the one carried out in the $n=2$ case we arrive at the Pfaf\/f\/ian system $\O_\mathfrak{F}$ (see \cite{Favard57} for details) whose integral curves are Frenet lifts of curves in $G/SL(3,\mathbb{R})$
\begin{gather*}
\O_\mathfrak{F}: \ \o^2=0,\ \ \o^3=0,\ \ \o^3_1=0, \ \ \o^2_1-\o^1=0,\ \ \o^3_2-\o^1=0,\ \ \o^1_1=0,\nonumber\\
\phantom{\O_\mathfrak{F}:} \  \ \o^2_2=0,\ \ \o^2_3-3\o^1_2=0.
\end{gather*}
I want to show that a Cartan prolongation of $\O_\mathfrak{F}$ is a contact system. We calculate that
\begin{equation*}
\O_\mathfrak{F}^\perp=\big\{v_1+v_5+v_9, \,v_7+3v_{11}, \, v_{10}\big\}.
\end{equation*}
The Cartan prolongation we shall consider is obtained from this:
\begin{equation*}
\widehat{\O}_\mathfrak{F}^\perp: \ \big\{v_1+v_5+v_9+\k_1(v_7+3v_{11})+\k_2v_{10},\, \P {\k_1},\, \P {\k_2}\big\}
\end{equation*}
def\/ined over $G\times\mathbb{R}^2$. To apply the generalised Goursat normal form, Theorem \ref{genGoursatNF_uniform}, we work with this dual bundle $\CV:=\widehat{\O}_\mathfrak{F}^\perp$ and
calculate
\begin{center}
\begin{tabular}{c|ll}
$i$ & $\CV^{(i)}$ & $\ch\CV^{(i)}$\cr
\cline{1-3}
0&$\big\{v_1+v_5+v_9+\k_1(v_7+3v_{11})+\k_2v_{10},\ \P {\k_1},\P {\k_2}\big\}$&$\{0\}$\cr
1&$\CV\oplus\big\{v_7+3v_{11},\ v_{10}\big\}$   & $\{\P {\k_1},\ \P {\k_2}\}$\cr
2&$\CV^{(1)}\oplus\big\{v_4+2v_8,\ v_7-v_{11}\big\}$ & $\ch\CV^{(1)}\oplus\{v_7+3v_{11},\ v_{10}\}$\cr
3&$\CV^{(2)}\oplus\big\{v_4-2v_8,\ v_1+v_5-5v_9\big\}$ & $\ch\CV^{(2)}\oplus\{v_4+2v_8,\ v_7-v_{11}\}$\cr
4&$\CV^{(3)}\oplus\big\{v_6,\ v_1-3v_5\big\}$ &  $\ch\CV^{(3)}\oplus\{v_4-2v_8,\ v_9\}$\cr
5&$T(G\times\mathbb{R}^2)$ & $T(G\times\mathbb{R}^2)$\cr
\end{tabular}
\end{center}

So hypotheses a), b) and c) of Theorem \ref{genGoursatNF_uniform} are satisf\/ied with $q=2$ and derived length $k=5$. Since $q>1$, it remains to check the singular variety of the quotient $\CV^{(4)}/\ch\CV^{(4)}$. From the table we see that
\[
\ch\CV^{(4)}=\big\{v_4,\ v_7,\ v_8,\ v_9,\ v_{10},\ v_{11},\ \P {\k_1},\ \P {\k_2}\big\}
\]
and
\[
\CV^{(4)}=\big\{v_1,\ v_4,\ v_5,\ v_6,\ v_7,\ v_8,\ v_9,\ v_{10},\ v_{11},\ \P {\k_1},\ \P {\k_2}\big\}
\]
and hence
\[
\widehat{\CV}^{(4)}:=\CV^{(4)}\big/\ch\CV^{(4)}=\big\{[v_1],\ [v_5],\ [v_6]\big\}.
\]
We obtain that $\text{Im}\,\hat{\delta}_4=\{[v_2],\ [v_3]\}$ and the polar matrix of the point $\langle a_1[v_1]+a_2[v_5]+a_3[v_6]\rangle\in\mathbb{P}\widehat{\CV}^{(4)}$ is
\[
\left(\begin{matrix}
-a_2&a_1&0\cr
-a_3&0&a_1\cr
\end{matrix}\right)
\]
which has unit rank if and only if $a_1=0$. Hence the singular variety of $\widehat{\CV}^{(4)}$ is $\mathbb{R}\mathbb{P}^1$ with singular bundle
$\widehat{\mathcal{B}}=\big\{[v_5],  [v_6]\big\}$. Consequently, the resolvent bundle in this case is
\[
\mathcal{R}\left(\CV^{(4)}\right)=\{v_4,\ v_5,\ v_6,\ v_7,\ v_8,\ v_9,\ v_{10},\ v_{11}\}\oplus\{\P {\k_1},\ \P {\k_2}\}=\mathfrak{s}\mathfrak{l}(3,\mathbb{R})\oplus\mathbb{R}^2.
\]
The resolvent bundle is integrable, showing that $\CV^{(4)}$ admits an integrable Weber structure, fulf\/illing hypothesis d) of Theorem~\ref{genGoursatNF_uniform}. We can therefore conclude that $\widehat{\O}_\mathfrak{F}^\perp$ is locally equivalent to the contact distribution $\mathcal{C}^{(5)}_2$ on jet space $J^5(\mathbb{R},\mathbb{R}^2)$.

For future reference, we note that we often abuse the term ``derived type" by referring to the list of lists
\[
\big[\,\big[\dim\mathcal{V}^{(j)},\, \dim\ch\CV^{(j)}\big]\,\big]_{j=0}^k
\]
as the {\it derived type} of $\CV$, where $k$ is its the derived length. Thus for the example just treated we may say that its derived type is
\[
\big[\,[3,0],[5,2],[7,4],[9,6],[11,8],[13,13]\,\big]
\]
since it is really the dimensions of these bundles that settles the {\it recognition} problem.

\subsection[Differential invariants via Contact]{Dif\/ferential invariants via {\itshape Contact}}\label{section4.1}

We now apply procedure {\it Contact} to construct the two equi-af\/f\/ine invariants for space curves. In the previous subsection we demonstrated that $\widehat{\O}_\mathfrak{F}^\perp$ is a uniform Goursat bundle. The f\/irst step in procedure {\it Contact} requires the $q+1=2+1=3$ invariants of the (integrable) resolvent bundle,~$\mathcal{R}(\CV^{(4)})$. The vector f\/ields spanning this bundle are all vertical for the projection $G\times\mathbb{R}^3\to \mathbb{R}^3$ and, indeed, they frame its f\/ibres. It follows that the functions annihilated by these vector f\/ields are spanned by $x$, $y$, $z$ -- the coordinates of the homogeneous space $G/SL(3,\mathbb{R})$ in which the curves are immersed. Note that no integration is required\footnote{See Section~\ref{section5} for further comment on this aspect.}.  A local coordinate calculation verif\/ies this claim.
By step~b) we are at liberty to take any one of these as the parameter along the curve. We take $x$ for this purpose. Continuing to follow~b) we f\/ix a vector f\/ield $\boldsymbol{Z}\in\CV$ such that $\boldsymbol{Z}x=1$. At this point we must construct the vector f\/ields that span $\CV$. This is straightforward since $\CV$ is constructed from the left-invariant vector f\/ields on~$G$ (as well as~$\P {\k_1}$,~$\P {\k_2}$). We obtain
\begin{gather*}
 v_1=a_1\P x+a_4\P y+a_7\P z,\qquad v_2=a_2\P x+a_5\P y+a_8\P z,
\qquad v_3=a_3\P x+a_6\P y+a_9\P z,\\
 v_4=a_1\P {a_1}-a_3\P {a_3}+a_4\P {a_4}-a_6\P {a_6}+a_7\P {a_7},\qquad
v_5=a_2\P {a_1}+a_5\P {a_4}+a_8\P {a_7},\\
 v_6=a_3\P {a_1}+a_6\P {a_4}+a_9\P {a_7},\qquad
v_7=a_1\P {a_2}+a_4\P {a_5}+a_7\P {a_8},\\
 v_8=a_2\P {a_2}-a_3\P {a_3}+a_5\P {a_5}-a_6\P {a_6}+a_8\P {a_8},\qquad
v_9=a_3\P {a_2}+a_6\P {a_5}+a_9\P {a_8},\\
 v_{10}=a_1\P {a_3}+a_4\P {a_6},\qquad v_{11}=a_2\P {a_3}+a_5\P {a_6},\qquad
v_{12}=\P {\k_1},\qquad v_{13}=\P {\k_2},
\end{gather*}
where
\[
a_9=\frac{1-a_2a_6a_7+a_1a_6a_8-a_3a_4a_8+a_3a_5a_7}{a_1a_5-a_4a_2}.
\]
From these we obtain
\[
\boldsymbol{Z}=\frac{1}{a_1}\big(v_1+v_5+v_9+\k_1(v_7+3v_{11})+\k_2v_{10}\big)\in\CV.
\]
Note that $\boldsymbol{Z}$ is the total dif\/ferential operator in this example. Finally, we let $z^1_0=y$, $z^2_0=z$ as in step~c) and compute the higher order coordinates as in step~d) by dif\/ferentiation by $\boldsymbol{Z}$. This construction provides the components of the local equivalence $\phi:G\times\mathbb{R}^2\to J^5(\mathbb{R},\mathbb{R}^2)$ that identif\/ies $\CV$ with the contact distribution $\mathcal{C}^{(5)}_2$ on $J^5(\mathbb{R},\mathbb{R}^2)$. Local inverse $\psi$ of the map $\phi$ gives the explicit dif\/ferential invariants $\k_1$, $\k_2$, for curves in $\mathbb{R}^3$ up to  ``equi-af\/f\/ine motions'':
\begin{gather*}
\k_1=\Big( 24z_2^2w_3w_5-35z_2^2w_4^2-60z_2w_3^2z_4+60z_2w_3z_3w_4-24z_2w_2w_3z_5+70z_2w_4w_2z_4\\
\phantom{\k_1=}{}  -24z_2w_2z_3w_5+60w_2w_3z_3z_4-60w_2z_3^2w_4- 35w_2^2z_4^2+24z_5w_2^2z_3\Big)\big/(w_3z_2-w_2z_3)^{7/3},
\\
\k_2=\Big( -18z_2^3w_3w_4w_5+25z_2^3w_4^3-36z_2^2w_3^3z_5+90z_2^2w_3^2z_4w_4+36z_2^2w_3^2z_3w_5\\
\phantom{\k_2=}{} -90z_2^2w_3z_3w_4^2+18z_2^2w_3w_2z_4w_5+18z_2^2w_3w_2w_4z_5-
  75z_2^2w_2w_4^2z_4+18z_2^2w_2z_3w_4w_5\\
\phantom{\k_2=}{}  -90z_2w_3^2w_2z_4^2+72z_2w_3^2z_5w_2z_3-
 18z_2w_3w_2^2z_4z_5-72z_2w_3w_2z_3^2w_5+90z_2w_2z_3^2w_4^2\\
\phantom{\k_2=}{} +75z_2w_2^2w_4z_4^2-
  18z_2w_2^2z_3w_4z_5-18z_2w_2^2z_3z_4w_5+90w_3w_2^2z_3z_4^2-36w_3z_5w_2^2z_3^2\\
\phantom{\k_2=}{}-25w_2^3z_4^3-90w_2^2z_3^2w_4z_4+18w_2^3z_3z_4z_5+36w_2^2z_3^3w_5\Big)\big/(w_3z_2-w_2z_3)^{7/2},
\end{gather*}
where $z_i=f^{(i)}(x)$, $w_j=g^{(j)}(x)$ for arbitrary smooth functions $f(x)$, $g(x)$ such that
\[g{'''}(x)f{''}(x)-g{''}(x)f{'''}(x)\neq 0.\]
Thus, if the curve in question has parametrisation $(x,f(x),g(x))$ then the complete set of dif\/ferential invariants is $\{\k_1,\k_2\}$. Pulling back the remaining nonzero elements of the Maurer--Cartan form on $G$ by $\psi$ give explicit expressions for invariant dif\/ferential one-forms on $J^{\,5}(\mathbb{R},\mathbb{R}^2)$. One can easily obtain the invariants in an arbitrary parametrisation of the curve or in terms of arclength parametrisation from these formulas but we won't record these here\footnote{Note that the inversion of $\phi$ is greatly assisted by the fact that the algebraic equations involved are guaranteed, by the proof of Theorem~3.3 in~\cite{Vassiliou1} to be block triangular.}.

\section{Curves in Riemannian manifolds} \label{section5}

The geometries discussed in previous sections of this paper and the one treated in~\cite{Vassiliou2} are all of {\it Klein} type and one may wonder how the method fares when curvature is introduced. Also, in our previous (illustrative) examples, we have permitted ourselves knowledge of the explicit transitive group action in order to construct the dif\/ferential system $\CV$ and we do {\it not} want to make this assumption in general. In this paper we will be content to establish a result for Riemannian geometry and illustrate this by two examples. However, we conjecture that similar results hold for other Cartan geometries.

Let $\o^i$, $\pi^i_j$ be the components of the Cartan connection\footnote{See Sharpe~\cite{Sharpe} for an account of Klein--Cartan geometries and Cartan connections; see also~\cite[Chapter~7]{Spivak79}.} of an arbitrary Riemannian mani\-fold~$(M,g)$, where
$\dim M=n$ and
$\o^i$, $1\leq i\leq n$ are semi-basic 1-forms for the projection $\mathcal{F}(M)\to M$;  $\mathcal{F}(M)$ is the orthonormal frame bundle over~$M$.
The structure equations for the coframe on~$\mathcal{F}(M)$ are
\begin{equation*}
d\o^i=\sum_{j=1}^n\,\o^j\wedge\pi^i_j,\qquad
 d\pi^i_j-\sum_{k=1}^n\,\pi^i_k\wedge\pi^k_j=\frac{1}{2}
\sum_{k,l=1}^n\,R^{\,i}_{j\,k\,l}\,\o^k\wedge\o^l,
\end{equation*}
where all indices range from 1 to $n$ and $\pi^i_j+\pi^j_i=0$ for all $i$, $j$.
Dually, the frame satisf\/ies
\begin{gather}\label{StrEqSO5}
[\P {\o^i},\P {\o^j}]=-\frac{1}{2}\sum_{k,l=1}^nR^k_{lij}\P {\pi^k_l},\qquad [\P {\o^i},\P {\pi^l_j}]=-\delta_{ij}\P {\o^l},\qquad
 [\P {\pi^i_l},\P {\pi^j_k}]=\delta_{kl}\P {\pi^i_j}+\delta_{ij}\P {\pi^l_k}.\!\!\!
\end{gather}

We now state the main result of this section.
\begin{theorem}\label{curvesRiemannian}
Let $n\geq 3$ and $\g:I\to M$ be an immersed curve in an $n$-dimensional Riemannian manifold $(M,g)$, $I\subseteq\mathbb{R}.$ Let
$\o^i,\pi^i_j$ be the components of the Cartan connection for $(M,g)$ on the orthonormal frame bundle $\mathcal{F}(M)\to M$.  Then
\begin{enumerate}\itemsep=0pt
\item[{\rm 1.}] There is a unique integral submanifold of
\begin{gather}
 \mathcal{V}=\Bigg\{\P {\o^1}+\kappa^1_0\P {\pi^2_1}+
\kappa^2_0\P {\pi^3_2}+\cdots+\kappa^{n-1}_0\P {\pi^n_{n-1}}+
\sum_{l=1}^{n-2}\kappa^1_l\P {\kappa^1_{l-1}}\nonumber\\
\hphantom{\mathcal{V}=\Bigg\{}{} +\sum_{l=1}^{n-3}\kappa^2_l\P {\kappa_{l-1}}+\sum_{l=1}^{n-4}\kappa^3_l\P {\kappa_{l-1}}+\cdots+
\kappa^{n-2}_1\P {\kappa^{n-2}_0},\ \P {\kappa^1_{n-2}},
\P {\kappa^2_{n-3}},\ldots, \P {\kappa^{n-1}_0}  \Bigg\},\!\!\!\label{RiemCurvesBundle}
\end{gather}
which projects to the Frenet lift of $\g$ to the orthonormal frame bundle $\mathcal{F}(M)$.

\item[{\rm 2.}] The sub-bundle $\CV$ is a uniform Goursat bundle with signature
\[
\langle\, \stackrel{n-1}{\overbrace{0,\, 0,\ldots,0}},\, n-1\,\rangle.
\]
Hence there is a local diffeomorphism $\phi$ which identifies it with the contact distribution~$\mathcal{C}^{(n)}_{n-1}$ on jet space $J^n(\mathbb{R},\mathbb{R}^{n-1})$.

\item[{\rm 3.}] The local diffeomorphism $\phi$ can be constructed by differentiation and algebraic operations alone.

\item[{\rm 4.}] If the isometries of $(M,g)$ act transitively, then the local diffeomorphism $\phi$ induces local coordinate formulas for the complete invariants $\kappa^1_0,\kappa^2_0,\ldots,\kappa^{n-1}_0$ of the curve $\g$ up to isometries of $(M,g)$.
\end{enumerate}
\end{theorem}

\begin{proof} The proofs of~1 and~4 are similar to the proof of Theorem~\ref{2dequiaffineCongruence}. The proof of~2 is complicated to write down in all generality; it is more enlightening to prove it in a suf\/f\/iciently non-trivial case. We therefore write out the detail proof in the case $n=6$, after which the general case will be clear.

The distribution in question, for $n=6$, is
\begin{gather*}
 \mathcal{V}=\Bigg\{\P {\o^1}+\sum_{l=1}^{5}\kappa^l_0\P {\pi^l_{l+1}}+\sum_{l=1}^4\kappa^1_l\P {\kappa^1_{l-1}}+\sum_{l=1}^3\kappa^2_l\P {\kappa^2_{l-1}}+
\sum_{l=1}^2\kappa^3_l\P {\kappa^3_{l-1}}+\kappa^4_1\P {\kappa^4_0},\nonumber\\
\hphantom{\mathcal{V}=\Bigg\{}{} \P {\kappa^1_4},\ \P {\kappa^2_3}, \ \P {\kappa^3_2},\ \P {\kappa^4_1},\ \P {\kappa^5_0}\Bigg\}.
\end{gather*}
We are required to prove that $\mathcal{V}\simeq\mathcal{C}^{(6)}_{5}$. Let us arrange some of the frame elements into subsets as follows
\begin{gather*}
 \mathfrak{p}_1:\ \big\{\P {\pi^5_6}\big\},\qquad
 \mathfrak{p}_2:\ \big\{\P {\pi^4_5},\ \P {\pi^4_6}\big\},\qquad
\mathfrak{p}_3:\ \big\{\P {\pi^3_4},\ \P {\pi^3_5},\ \P {\pi^3_6}\big\},\nonumber\\
\mathfrak{p}_4:\ \big\{\P {\pi^2_3},\ \P {\pi^2_4},\ \P {\pi^2_5},\ \P {\pi^2_6}\big\},\qquad
\mathfrak{p}_5:\ \big\{\P {\pi^1_2},\ \P {\pi^1_3},\ \P {\pi^1_4},\ \P {\pi^1_5},\ \P {\pi^1_6}\big\}.
\end{gather*}
It is easy to show from (\ref{StrEqSO5}) that
\[
\mathfrak{h}_\ell=\bigoplus_{j=1}^\ell\,\mathfrak{p}_j\subseteq \mathfrak{so}(6),\qquad 1\leq \ell\leq 5,
\]
determine a f\/lag of Lie subalgebras. Using this we compute that for each $s$ in the range $1\leq s\leq 6$, the quotients $\widehat{\mathcal{V}}^s:=\mathcal{V}^{(s)}/\mathcal{V}^{(s-1)}$ have basis representatives
\begin{gather*}
\widehat{\mathcal{V}}^{(1)}=\big\{\P {\k^1_3},\ \P {\k^2_2},\ \P {\k^3_1},\ \P {\k^4_0},\ \P {\pi^5_6}\big\},\qquad
\widehat{\mathcal{V}}^{(2)}=\big\{\P {\k^1_2},\ \P {\k^2_1},\ \P {\k^3_0},\ \P {\pi^4_5},\ \P {\pi^4_6}\big\},\nonumber\\
\widehat{\mathcal{V}}^{(3)}=\big\{\P {\k^1_1},\ \P {\k^2_0},\ \P {\pi^3_4},\ \P {\pi^3_5},\ \P {\pi^3_6}\big\},\qquad
\widehat{\mathcal{V}}^{(4)}=\big\{\P {\k^1_0},\ \P {\pi^2_3},\ \P {\pi^2_4},\ \P {\pi^2_5},\ \P {\pi^2_6}\big\},\nonumber\\
\widehat{\mathcal{V}}^{(5)}=\big\{\P {\pi^1_2},\ \P {\pi^1_3},\ \P {\pi^1_4},\ \P {\pi^1_5},\ \P {\pi^1_6}\big\},\qquad
\widehat{\mathcal{V}}^{(6)}=\big\{\P {\o^2},\,\P {\o^3},\,\P {\o^4},\,\P {\o^5},\ \P {\o^6}\big\}.
\end{gather*}
The Cauchy systems have bases
\begin{gather*}
\ch{\mathcal{V}}^{(1)}=\big\{\P {\k^1_4},\ \P {\k^2_3},\ \P {\k^3_2},\ \P {\k^4_1},\ \P {\k^5_0}\big\},\nonumber\\
\ch{\mathcal{V}}^{(2)}=\big\{\P {\k^1_3},\ \P {\k^2_2},\ \P {\k^3_1},\ \P {\k^4_0},\ \P {\pi^5_6}\big\}\oplus\ch\mathcal{V}^{(1)},\nonumber\\
\ch{\mathcal{V}}^{(3)}=\big\{\P {\k^1_2},\ \P {\k^2_1},\ \P {\k^3_0},\ \P {\pi^4_5},\ \P {\pi^4_6}\big\}\oplus\ch\mathcal{V}^{(2)},\nonumber\\
\ch{\mathcal{V}}^{(4)}=\big\{\P {\k^1_1},\ \P {\k^2_0},\ \P {\pi^3_4},\ \P {\pi^3_5},\ \P {\pi^3_6}\big\}\oplus\ch\mathcal{V}^{(3)},\nonumber\\
\ch{\mathcal{V}}^{(5)}=\big\{\P {\k^1_0},\ \P {\pi^2_3},\ \P {\pi^2_4},\ \P {\pi^2_5},\ \P {\pi^2_6}\big\}\oplus\ch\mathcal{V}^{(4)}
\end{gather*}
and hence
\[
\mathcal{V}^{(5)}/\ch\mathcal{V}^{(5)}=
\big\{\big[\P {\o^1}\big],\ \big[\P {\pi^1_2}\big],\ \big[\P {\pi^1_3}\big],\
\big[\P {\pi^1_4}\big],\ \big[\P {\pi^1_5}\big],\big[\ \P {\pi^1_6}\big]\big\}.
\]
Structure equations (\ref{StrEqSO5}) show that the singular sub-bundle is
\[
\widehat{\mathcal{B}}=\big\{\big[\P {\pi^1_2}\big],\ \big[\P {\pi^1_3}\big],\
\big[\P {\pi^1_4}\big],\ \big[\P {\pi^1_5}\big],\ \big[\P {\pi^1_6}\big]\big\}
\]
and therefore the resolvent bundle is given by
\[
\mathcal{R}\big(\CV^{(5)}\big)=\big\{\P {\pi^i_j}\big\}\oplus\big\{\P {\k^a_b}\big\}=\mathfrak{so}(6)\oplus\mathbb{R}^{15},
\]
where $a$, $b$, $i$, $j$ range over all possible values. These calculations show that that the derived type of $\CV$ is
\[
[[6,0],[11,5],[16,10],[21,15],[26,20],[31,25],[36,36]]
\]
from which one deduce's that the signature of $\CV$ is
$
\langle 0, 0, 0, 0, 0, 5\rangle.
$
Since $\mathcal{R}(\CV^{(5)})$ is integrable and all other hypotheses of Theorem \ref{genGoursatNF_uniform} are satisf\/ied with $q=5$ and $k=6$, we have shown that~$\mathcal{V}$ is locally equivalent to the contact system $\mathcal{C}^{(6)}_5$ on $J^6(\mathbb{R},\mathbb{R}^5)$.

To prove 3, we invoke Theorem~4.2 of \cite{Vassiliou2} which shows that to construct the local equivalence~$\phi$ identifying $\mathcal{V}$ with $\mathcal{C}^{(6)}_5$ one only requires a complete set of invariants of the resolvent bundle~$\mathcal{R}(\CV^{(5)})$. However, it is elementary to see that no element of $\mathcal{R}(\CV^{(5)})$ has components tangent to~$M$; on the other hand $\mathcal{R}(\CV^{(5)})$ spans the tangent spaces of the f\/ibres over~$M$. Hence, any coordinate system on $M$ provides the needed invariants and no integration need be performed.
\end{proof}

In case the Riemannian manifold $M$ does {\it not have} a transitive isometry group then $\mathcal{V}$ may simply be regarded as a {\it control system} on $M$ with controls $\kappa^1_{n-2}, \kappa^2_{n-3},\allowbreak\ldots,\kappa^{n-1}_0$ and with all other coordinates on $E$ playing the role of state variables ({\it outputs}, in the language of control theory). The arc-length parametrisation along the curve plays the role of {\it time} in this control theoretic interpretation. The fact that $\mathcal{V}$ is a Goursat bundle then proves that it is {\it differentially flat}, a type of control system currently under investigation in the control community. In fact it is currently a signif\/icant open problem in control theory to geometrically characterise all dif\/ferentially f\/lat control systems. One other way to state Theorem~\ref{curvesRiemannian} is to say that the ``natural" framing of curves in any Riemannian manifold is dif\/ferentially f\/lat.

If $M$ does possess a transitive isometry group $G$ then the functions $\kappa^1_0,\ldots,\kappa^{n-1}_0$ form a~complete set of curve invariants up to the action of~$G$.

One may wonder whether Theorem \ref{curvesRiemannian} is a genuine advance in the theory and/or practice of moving frames either in the sense of Cartan or in the sense of Fels--Olver. The contention of this paper is that it does represent an advance in both the theory and the practice of moving frames. As to the theory, I argue that endowing curves in a Cartan geometry with an explicit contact structure is signif\/icant given the fundamental role that contact structures play in geometry and dif\/ferential equations. As to the practice of the method of moving frames, in the case of Riemannian manifolds $(M,g)$ we have given a framing for curves in~$M$ as solutions of a~dif\/ferential system $\widehat{\O}_{\mathfrak{F}}$ for any metric $g$ and shown that it can be explicitly identif\/ied with the contact system on a jet space $J^k(\mathbb{R},\mathbb{R}^q)$ for some~$k$,~$q$. In case the isometry group of $(M,g)$ acts transitively, this identif\/ication delivers the complete set of invariant data including the dif\/ferential curve invariants. We formulate this construction as an algorithm.

\bigskip

\centerline{\it Algorithm Riemannian curves}
\begin{enumerate}\itemsep=0pt
\item[] $\mathbf{INPUT:}$ Riemannian manifold, $(M,g)$, dimension $n\geq 3$.
\item[a)] Construct a parametrisation of $SO(n)$  and orthonormal coframe $\o$ on $M$.
\item[b)] Lift $\o$ to the orthonormal frame bundle over $M$  and build the Cartan connection $\Theta$ for $(M,g)$.
This involves linear algebra.
\item[c)] From $\Theta$, build the dif\/ferential system $\mathcal{V}$ def\/ined  in Theorem \ref{curvesRiemannian}.
\item[d)] Apply procedure {\it Contact} using, as invariants of  the resolvent bundle, any coordinate system on $M$.
\item[e)] Construct dif\/feomorphism $\phi$ by procedure \textit{Contact A}.
\item[] $\mathbf{OUTPUT:}$ 
Local dif\/feomorphism $\phi$ identifying $\CV$ with  contact distribution $\mathcal{C}^{(n)}_{n-1}$.
In case the isometry  group of $(M,g)$ acts  transitively,  obtain  complete,  explicit  invariant data for curves in $M$, following the inversion of $\phi$.
\end{enumerate}

Note that {\it Riemannian curves} is an algorithm in as much as steps a)--e) do not require {\it any integration} to be performed. All the steps involved are algebraic.

\subsection{Curves in the Poincar\'e half-space}\label{section5.1}

The aim of this subsection is to apply the previous theorem to construct explicit expressions for curvature and torsion for curves in the Poincar\'e half-space $H^3$, with Riemannian metric
\[
g=\frac{dx^2+dy^2+dz^2}{z^2}.
\]
In principle we could approach this by putting coordinates on the Lie group $SO(3,1)$ and then consider curves in the homogeneous space $SO(3,1)/SO(3)$, as we did in the case of equi-af\/f\/ine space curves. However, this leads to unwieldy expressions which are dif\/f\/icult to handle, even with the help of a computer. Instead, we will use the method of equivalence to construct the Cartan connection for $g$ and then use this to build the canonical Pfaf\/f\/ian system $\O_\mathfrak{F}$ for the Frenet frame of a generic curve in $H^3$. In fact, we will construct the dual vector f\/ield distribution
$\mathcal{V}=\O_\mathfrak{F}^\perp$.

Begin by lifting the 1-forms
\[
\theta^1=\frac{dx}{z},\qquad \theta^2=\frac{dy}{z},\qquad \theta^3=\frac{dz}{z}
\]
to the orthonormal frame bundle $\mathcal{F}(H^3)$ over $H^3$ by
\[
\left(\begin{matrix} \o^1\cr \o^2\cr \o^3\end{matrix}\right)=O\left(\begin{matrix}
\theta^1\cr\theta^2\cr \theta^3
\end{matrix}\right),
\]
where
\[
O=\left(\begin{matrix} \cos b \cos c-\sin a \sin b \sin c & -\cos a \sin c & \sin a \cos b \sin c+
\sin b \cos c\\
\cos b \sin c+\sin a \sin b \cos c & \cos a \cos c &\sin b \sin c-\sin a \cos b \cos c\\
-\cos a \sin b &\sin a&\cos a \cos b
\end{matrix}\right)
\]
parametrises $SO(3)$. The f\/ibres of the orthonormal frame bundle $\mathcal{F}(H^3)\to H^3$ are dif\/feomorphic to $SO(3)$ whose Maurer--Cartan form is the $\mathfrak{s}\mathfrak{o}(3)$-valued 1-form
\[
\O=\left(\begin{matrix}
                    0 & \o^1_2 & \o^1_3\\
                    \o^2_1 & 0 & \o^2_3\\
                    \o^3_1 & \o^3_2 & 0\end{matrix}\right),
\]
where
\begin{gather*}
\o^1_2=-dc-\sin a\;db,\qquad
\o^1_3=\sin c\;da+\cos a \cos c\;db,\qquad
\o^2_3=\sin c\,\cos a\;db-\cos c\;da,
\end{gather*}
and $\o^j_i+\o^i_j=0$. As is usual in the method of equivalence, we compute the structure equations of the semi-basic forms  obtaining
\begin{gather*}
 d\o^1=\o^1_2\wedge\o^2+\o^1_3\wedge\o^3+(\sin c \sin b-\cos c\sin a\cos b) \o^1\wedge\o^2+
 \cos a\cos b\;\o^1\wedge\o^3,\\
d\o^2=\o^2_1\wedge\o^1+\o^2_3\wedge\o^3-(\cos c \sin b+\sin a \cos b \sin c) \o^1\wedge\o^2+
 \cos a \cos b\;\o^2\wedge\o^3,\\
d\o^3=\o^3_1\wedge\o^1+\o^3_2\wedge\o^2-(\sin b \cos c+\sin a \cos b \sin c) \o^1\wedge\o^3\\
\phantom{d\o^3=}{} - (\sin b \sin c-\sin a \cos b \cos c) \o^2\wedge\o^3.
\end{gather*}
All torsion can be absorbed by redef\/ining connection forms
\begin{gather*}
 \pi^1_2=\o^1_2+(\sin c \sin b-\sin a \cos b \cos c) \o^1-
(\sin a \cos b \sin c+\sin b\cos c) \o^2,\\
 \pi^1_3=\o^1_3+\cos a \cos b \o^1-(\sin a \cos b \sin c+\sin b \cos c   ) \o^3,\\
 \pi^2_3=\o^2_3+\cos a \cos b \o^2-(\sin b \sin c-\sin a \cos b \cos c) \o^3,
\end{gather*}
and we obtain structure equations
\begin{gather*}
 d\o^1= \pi^1_2\wedge \o^2+\pi^1_3\wedge\o^3,\qquad
 d\o^2=-\pi^1_2\wedge\o^1+\pi^2_3\wedge\o^3,\qquad
 d\o^3=-\pi^1_3\wedge\o^1-\pi^2_3\wedge\o^2.
\end{gather*}
The remaining structure equations are
\begin{gather*}
d\pi^1_2=-\pi^1_3\wedge\pi^2_3+\o^1\wedge\o^2,\qquad
d\pi^1_3= \pi^1_2\wedge\pi^2_3+\o^1\wedge\o^3,\qquad
d\pi^2_3=-\pi^1_2\wedge\pi^1_3+\o^2\wedge\o^3.
\end{gather*}
These structure equations are those of the Maurer--Cartan form on $SO(3,1)$. Moreover, we are now able to def\/ine the $\mathfrak{e}(3)$-valued Cartan connection
\[
\widehat{\O}=\left(\begin{matrix} 0 & 0 & 0 & 0\\
                    \o^1 & 0 & \pi^1_2 & \pi^1_3\\
                    \o^2 & \pi^2_1 & 0 & \pi^2_3\\
                    \o^3 &  \pi^3_1 & \pi^3_2 & 0\end{matrix}\right),
\]
with curvature
\[
d\widehat{\O}+\widehat{\O}\wedge\widehat{\O}=
\left(\begin{matrix} 0 & 0 & 0 & 0\\
                     0 & 0 & \o^1\wedge \o^2 &  \o^1\wedge\o^3\\
                     0 & \o^2\wedge\o^1 & 0 &  \o^2\wedge\o^3\\
                     0 & \o^3\wedge\o^1 & \o^3\wedge\o^2 & 0
                     \end{matrix}\right)
\]
for metric $g$.

Setting $n=3$ in (\ref{RiemCurvesBundle}), and adopting the usual notation $\kappa^1_0=\kappa$, $\kappa^2_0=\tau$, $\kappa^1_1=\kappa_1$, we study the integral submanifolds of
\[
\mathcal{V}=\big\{\P {\o^1}+\kappa\P {\pi^1_2}+\tau\P {\pi^2_3}+\kappa_1\P {\kappa},\ \P {\kappa_1},\
\P {\tau}\big\}.
\]
It is easy to check that $\mathcal{V}$ has derived type
\[
[[3,0],[5,2],[7,4],[9,9]].
\]
It can be checked that $\rho_1=\rho_2=0$ and hence the signature of $\CV$ is $\langle 0, 0, 2\rangle$. This is the signature of contact system $\mathcal{C}^{(3)}_2$. To complete the check that $\mathcal{V}$ is dif\/feomorphic to $\mathcal{C}^{(3)}_2$, we compute the resolvent bundle and check its integrability.

We f\/ind that $\ch\mathcal{V}^{(2)}=\{\P {\pi^2_3},\P {\kappa},\P {\kappa_1},\P {\tau}\}$, and
\[
\mathcal{V}^{(2)}/\ch\mathcal{V}^{(2)}=\big\{\big[\P {\o^1}\big],
\ \big[\P {\pi^1_2}\big],\ \big[\P {\pi^1_3}\big]\big\},
\]
whose structure is
\[
[\P {\o^1},\P {\pi^1_2}]\equiv \P {\o^2},\quad [\P {\o^1},\P {\pi^1_3}]\equiv\P {\o^3},\quad
[\P {\pi^1_2},\P {\pi^1_3}]\equiv 0\mod\ch\mathcal{V}^{(2)}.
\]
Hence, the singular bundle is $\widehat{\mathcal{B}}=\{[\P {\pi^1_2}],[\P {\pi^1_3}]\}$ and consequently the resolvent bundle is
\[
\mathcal{R}(\CV^{(2)})=\big\{\P {\pi^1_2},\ \P {\pi^1_3},\ \P {\pi^2_3},\ \P {\kappa},\ \P {\kappa_1},\
\P {\tau}\big\}=\mathfrak{so}(3)\oplus\mathbb{R}^3,
\]
which is clearly integrable. An easy calculation in local coordinates verif\/ies that the invariants of~$\mathcal{R}(\CV^{(2)})$ are indeed the coordinates $x$, $y$, $z$ on the base of the orthonormal frame bundle $\mathcal{F}(H^3)\to H^3.$

According to the Theorem we take one of $x$, $y$, $z$ as the independent variable. If  $y$ is taken for this purpose, then we can take $\mathcal{D}=\P {\o^1}+\kappa\P {\pi^1_2}+\tau\P {\pi^2_3}+\kappa_1\P {\kappa}$ as the direction of the total dif\/ferential operator since
$0\neq\mathcal{D}\;(y)=-z\cos a \sin c=\mu$. Then the operator of total dif\/ferentiation is
\[
\boldsymbol{Z}=\mu^{-1}\big(\P {\o^1}+\kappa\P {\pi^1_2}+\tau\P {\pi^2_3}+\kappa_1\P {\kappa}\big).
\]
Setting $u=x,v=z$, the remaining contact coordinates for $\mathcal{V}$ are then provided by dif\/ferentiation
\[
u_1=\boldsymbol{Z}\,u,\quad v_1=\boldsymbol{Z}\,v,\quad u_2=\boldsymbol{Z}\,u_1,\quad v_2=\boldsymbol{Z}\,v_1,
\quad u_3=\boldsymbol{Z}\,u_2,\quad v_3=\boldsymbol{Z}\,v_2.
\]
The map $\phi$
\[
(x,y,z,a,b,c,\kappa,\kappa_1,\tau)\mapsto
(y,u,v,u_1,v_1,u_2,v_2,u_3,v_3)
\]
pushes $\mathcal{V}$ forward to $\mathcal{C}^{(3)}_2$, which, in contact coordinates has the form
\[
\big\{\P y+u_1\P u+v_1\P v+u_2\P {u_1}+v_2\P {v_1}+u_3\P {u_2}+v_3\P {v_2},\ \P {u_3},\ \P {v_3}\big\}.
\]
A local inverse of $\phi$ is easily constructed and provides all the invariant data for curves in~$H^3$. In particular we deduce explicit expressions for curvature and torsion for curves
$\gamma(y){=}(u(y), y, v(y)),\!$ where for instance $u_1$, $v_2$, denote derivatives $du/dy$, $d^2v/dy^2$, etc. We obtain
\begin{gather*}
 \kappa=-\Big(u_1^6+2u_1^4v_1^2+2u_1^4vv_2+3u_1^4-2u_1^3vu_2v_1+v_1^4u_1^2+2v_1^2u_1^2vv_2+4u_1^2vv_2+u_1^2v^2v_2^2\\
\phantom{\kappa=}{} +4v_1^2u_1^2+3u_1^2-2u_1v^2u_2v_1v_2-2u_1vu_2v_1^3-2u_1v_1u_2v+1+v_1^4+v^2u_2^2v_1^2+2v_1^2vv_2\\
\phantom{\kappa=}{} +2vv_2+v_2^2v^2+v^2u_2^2+2v_1^2\Big)^{1/2}\Big/\big(v_1^2+u_2^2+1\big)^{3/2},\\
\tau=\Big(3u_2^2u_1+3u_2v_1v_2+vu_2v_3-u_3u_1^2-u_3v_1^2-u_3vv_2-u_3\Big)v^2\Big/\Big(u_2^2v^2v_1^2+u_2^2v^2\\
\phantom{\tau=}{} -2u_2v^2v_1u_1v_2-2u_2vv_1u_1-2v_1^3u_1u_2v-2u_1^3v_1u_2v+u_1^2v^2v_2^2
+v^2v_2^2+4u_1^2vv_2\\
\phantom{\tau=}{}+2v_1^2vv_2+2vv_2+2u_1^4vv_2+2v_1^2u_1^2vv_2+1+2v_1^2+v_1^4+4u_1^2v_1^2+
3u_1^2+2u_1^4v_1^2\\
\phantom{\tau=}{}+v_1^4u_1^2+3u_1^4+u_1^6\Big).
\end{gather*}

Once again, the algebraic system that is presented for solution in this task is guaranteed to have a block triangular structure.

We remark that semi-circular arcs parallel to the $y-z$ plane
\[
\gamma(y)=\Big(C_1,\ y,\ \sqrt{C_2^2-y^2}\,\Big),\qquad C_2>0,\qquad -C_2<y<C_2
\]
are geodesics in $H^3$; accordingly it can be checked that $\k (\g(y) )=\t (\g(y) )\equiv 0$.

\subsection{Curves in constant curvature Riemannian 3-manifolds}\label{section5.2}

Let $\lambda$ be any nonzero real number. Here we point out that the construction of the previous subsection can be carried out for the three dimensional Riemannian manifold $(M(\lambda),g_\lambda)$ with metric
\[
g_\lambda=\frac{dx^2+dy^2+dz^2}{\left(1+\frac{\phantom{x}\lambda}{\phantom{x}4^{\phantom{2}}}(x^2+y^2+z^2)\right)^2}
\]
of constant curvature $\lambda$, where $M(\lambda)$ is an open subset of $\mathbb{R}^3$. Exactly the same calculation as before but with
\[
\theta^1=\frac{dx}{\L},\qquad \theta^2=\frac{dy}{\L},\qquad \theta^3=\frac{dz}{\L},
\]
where $\L=1+\frac{\phantom{x}\lambda}{\phantom{x}4^{\phantom{2}}}(x^2+y^2+z^2)$ gives rise to the corresponding Cartan connection
\[
\widehat{\O}=\left(\begin{matrix} 0 & 0 & 0 & 0\\
                    \o^1 & 0 & \pi^1_2 & \pi^1_3\\
                    \o^2 & \pi^2_1 & 0 & \pi^2_3\\
                    \o^3 &  \pi^3_1 & \pi^3_2 & 0\end{matrix}\right),
\]
with curvature
\[
d\widehat{\O}+\widehat{\O}\wedge\widehat{\O}=
-\lambda\left(\begin{matrix} 0 & 0 & 0 & 0\\
                     0 & 0 & \o^1\wedge \o^2 &  \o^1\wedge\o^3\\
                     0 & \o^2\wedge\o^1 & 0 &  \o^2\wedge\o^3\\
                     0 & \o^3\wedge\o^1 & \o^3\wedge\o^2 & 0
                     \end{matrix}\right).
\]
for the metric $g_\lambda$. Exactly the same calculation as the one carried out for the Poincar\'e half-space gives rise to the curvature and torsion for curves in $(M(\lambda),g_\lambda)$.  We get
\begin{gather*}
\kappa_\lambda=\lambda\frac{\sqrt{(yu_1-u)^2+(yv_1-v)^2+(vu_1-uv_1)^2}}{2\sqrt{u_1^2+v_1^2+1}},\\
\tau_\lambda=-\frac{\big(4+\lambda(y^2+v^2+u^2)\big)^2
\big(uv_3-vu_3+y(v_1u_3-v_3u_1)\big)}{8\lambda\big(1+u_1^2+v_1^2\big)\big( (yu_1-u)^2+(yv_1-v)^2+(vu_1-uv_1)^2  \big)}.
\end{gather*}
The point to note here as in the previous example is that we are not required to know the explicit formulas for the action of the isometries on $M(\lambda)$ before the invariants and the moving frame can be computed. Only inf\/initesimal data is required, in the form of the Cartan connection and then no integration need be performed.

\begin{remark}
Interestingly, $\tau_\lambda$ is not a continuous function of $\lambda$ at $\lambda=0$, while
$\lim\limits_{\lambda\to 0}\k_\lambda$ is not the curvature of the curve in the corresponding limiting metric
$\lim\limits_{\lambda\to 0}\,g_\lambda$, which is Euclidean.
\end{remark}



\section{Closing remarks}\label{section6}

In this paper we have demonstrated that curves in various geometries can be endowed with a contact geometry by combining Cartan's classical construction of moving frames with the generalised Goursat normal form. In particular we have shown that curves in {\it any} Riemannian manifold can be endowed with a contact geometry regardless of the nature of its isometry group or curvature tensor.

\looseness=1
We have also been concerned with the problem of explicitly computing dif\/ferential invariants of curves immersed in spaces equipped with a transitive action of a Lie group. If this Lie group action is explicitly known and not too complicated then the method of choice for computing the dif\/ferential invariants and other geometric data is the normalisation of the group action as in the equivariant moving frames method of Fels and Olver, described in \cite{FelsOlver98,FelsOlver99}. This method is very general, has a simple and elegant theoretical foundation and presents as simple a computational task as could be hoped for.
However, if the explicit group action is not known or it is known but too complicated to work with and if the goal is explicit expressions for dif\/ferential invariants then the Fels--Olver method can't readily be used to compute invariants explicitly\footnote{However, it is sometimes possible to simplify the task of normalising the group action by using the known action of a subgroup, a result due to Kogan,~\cite{Kogan2000}.}.
In this case, we have shown that for Riemannian manifolds $(M,g)$, the contact geometry can be fruitfully used to derive dif\/ferential invariants and this requires as input data only the metric~$g$ and a~rea\-lisation of the structure group $SO(n)$.  With this data the Cartan connection for $(M,g)$ can be constructed by linear algebra and dif\/ferentiation. Subsequently, Theorem~\ref{curvesRiemannian} provides an algorithm, {\it Riemannian curves}, for the curve invariants and, if required the Fels--Olver equivariant moving frame.

\looseness=1
The symbolic computational aspects of procedure {\it Contact} and algorithm {\it Riemannian curves} presented in this paper should be mentioned brief\/ly. The {\tt Maple} package  {\tt DifferentialGeometry} is ideally suited to the computation of all the relevant bundles and determining the derived type of any sub-bundle $\CV\subset TM$ over manifold $M$. For instance, the two Riemannian examples presented in Section~\ref{section5}, take only a few minutes to complete commencing only with the metric and realisation of matrix group $SO(3)$. Furthermore, the construction of dif\/ferential curve invariants requires the inverse of the local dif\/feomorphism $\phi$ produced by procedure {\it Contact}. The proof of correctness of {\it Contact} in \cite{Vassiliou2} shows that this algebraic problem will be block triangular. Thus the procedures discussed in this paper have quite good computational features.

However, of much greater signif\/icance stands the proposition that the Frenet frames along a curve and hence the curve itself can be endowed with a contact geometry. This {\it should} have signif\/icance not only for the geometry of curves but also for Cartan's method of moving frames as well as for the equivariant moving frames method of Fels and Olver. This is because contact systems are fundamental geometric objects and play a central role in dif\/ferential geometry and dif\/ferential equations. In fact, the construction of any contact system out of the components of a Cartan connection can in many ways {\it replace} or {\it complement} the step by step construction of moving frames championed by Cartan. What is more, a characterisation of contact systems in {\it arbitrary jet spaces} in the spirit of the Goursat normal form is known \cite{Bryant79,Yamaguchi82} and could be applied to study the geometry of submanifolds of dimension $p>1$ in general geometries as we have done here in the case $p=1$. This raises the interesting question of the extent to which the results of this paper can be extended to homogeneous spaces in general and how they are connected to existing theory such as \cite{Cartan37,Chern85,Favard57,Griffiths74,Jensen77,Green78,Sulanke79,AlekLychaginVino1991,Sharpe,FelsOlver98,FelsOlver99,Steltsova09}.

In this respect  it should be mentioned that the Fels--Olver theory of moving frames has application well beyond the explicit calculation of dif\/ferential invariants and moving frames. Much can be accomplished within the theory even without this explicit knowledge; see Mansf\/ield~\cite{Mansfield2010} for details. What we hope to have achieved in this paper is the presention of evidence sup\-porting the proposition that it is useful to enrich the philosophy and practice of moving frames by exploring its links with contact geometry.

Finally, we mention that an intriguing question is the relationship between the contact geo\-metry of curves as explained here and the integrable motion of curves in various ambient manifolds \cite{Ivey01,ChouQu02,ManKamp06,Olver09}.

\subsection*{Acknowledgements}

I am indebted to the anonymous referees for insightful comments and for corrections which greatly improved the paper. Any remaining errors are mine.

\pdfbookmark[1]{References}{ref}
\LastPageEnding


\begin{thebibliography}{99}

\footnotesize\itemsep=0pt

\bibitem{AlekLychaginVino1991}
 Alekseev D.V., Vinogradov A.M., Lychagin V.V.,
 Geometry~I, {\it Encycl. Math. Sci.}, Vol.~28, Spinger-Verlag, Berlin, 1991.

\bibitem{Bryant79}
 Bryant R.L.,
 Some aspects of the local and global theory of Pfaf\/f\/ian systems, PhD thesis, University of North Carolina, Chapel-Hill, 1979.

\bibitem{BC3G}
 Bryant R.L., Chern S.S., Gardner R.B., Goldschmidt H.L., Grif\/f\/iths P.A.,
 Exterior dif\/ferential systems,
{\it Mathematical Sciences Research Institute Publications}, Vol.~18, Springer-Verlag, New York, 1991.

\bibitem{Cartan35}
 Cartan \'E.,
 La m\'ethode du rep\`ere mobile, la th\'eorie des groupes continus et les espaces g\'en\'eralis\'es
Hermann \& Cie, Paris, 1935.

\bibitem{Cartan37}
 Cartan \'E.,
 La th\'eorie des groupes f\/inis et continus et la g\'eom\'etrie dif\/f\'erentielle traitees par la m\'ethode du rep\`ere mobile,  Gautier-Villars, Paris, 1937.

\bibitem{Chern85}
 Chern S.S.,
 Moving frames, in  The Mathematical Heritage of \'Elie Cartan (Lyon, 1984),
 {\it Ast\'erique}, Vol.~1985, Numero Hors Serie, 1985, 67--77.

\bibitem{ChouQu02}
 Chou K.-S., Qu C.-Z.,
 Integrable equations arising from motions of plane curves,
 {\it Phys. D}  {\bf 162}  (2002), 9--33.

\bibitem{Favard57}
Favard J.,
Cours de g\'eom\'etrie dif\/f\'erentielle locale, Gauthier-Villars, Paris, 1957.

\bibitem{FelsOlver98}
 Fels M., Olver P.J.,
 Moving coframes. I.~A practical algorithm,
 {\it Acta Appl. Math.} {\bf 51}  (1998), 161--213.

\bibitem{FelsOlver99}
 Fels M., Olver P.J.,
 Moving coframes. II.~Regularization and theoretical foundations,
 {\it Acta Appl. Math.} {\bf 55}  (1999), 127--208.

\bibitem{Green78}
 Green M.L.,
 The moving frame, dif\/ferential invariants and rigidity theorems for curves in homogeneous spaces, {\it Duke Math. J.}  {\bf 45} (1978), 735--779.

\bibitem{Griffiths74}
 Grif\/f\/iths P.A.,
 On Cartan's method of Lie groups and moving frames as applied to uniqueness and existence questions in dif\/ferential geometry,
 {\it Duke Math. J.} {\bf 41}  (1974), 775--814.

\bibitem{Ivey01}
Ivey T.,
Integrable geometric evolution equations for curves,
in The Geometrical Study of Dif\/ferential Equations (Washington, DC, 2000), {\it Contemp. Math.}, Vol.~285, Amer. Math. Soc., Providence, RI, 2001, 71--84.

\bibitem{Jensen77}
 Jensen G.R.,
 Higher order contact of submanifolds of homogeneous spaces, {\it Lecture Notes in Mathematics}, Vol.~610, Springer-Verlag, Berlin~-- New York, 1977.

\bibitem{Kogan2000}
 Kogan I.A.,
 Inductive approach to moving frames and applications in classical invariant theory, PhD Thesis, University of Minesota, 2000.

\bibitem{Mansfield2010}
 Mansf\/ield E.L.,
 A guide to the symbolic invariant calculus, Cambridge University Press, Cambridge, to appear.

\bibitem{ManKamp06}
 Mansf\/ield E.L., van der Kamp P.E.,
 Evolution of curve invariants and lifting integrability, {\it J. Geom. Phys.} {\bf 56} (2006), 1294--1325.

\bibitem{Olver01}
 Olver P.J.,
 Moving frames -- in geometry, algebra, computer vision and numerical analysis, in Foundations of Computational Mathematics (Oxford, 1999), Editors R.~DeVore, A.~Iserles and E.~S\"uli, {\it London Math. Soc. Lecture Note Ser.}, Vol.~284, Cambridge University Press,  Cambridge, 2001, 267--297.

\bibitem{Olver09}
 Olver P.J.,
 Invariant submanifold f\/lows,
{\it J. Phys. A: Math. Theor.} {\bf 41} (2008), 344017, 22~pages.

\bibitem{ShadSluis}
Shadwick W.F., Sluis W.M.,
Dynamic feedback for the classical geometries,
in Dif\/ferential Geometry and Mathematical Physics (Vancouver, BC, 1993), {\it Contemp. Math.}, Vol.~170, Amer. Math. Soc., Providence, RI, 1994, 207--213.

\bibitem{Sharpe}
 Sharpe R.,
Dif\/ferential geometry. Cartan's generalisation of Klein's erlangen program, {\it Graduate Texts in Mathematics}, Springer-Verlag, New York, 1997.

\bibitem{Spivak70}
Spivak M.,
 A comprehensive introduction to dif\/ferential geometry, Vol.~1, Publish or Perish Press, 1970.

\bibitem{Spivak79}
 Spivak M.,
A comprehensive introduction to dif\/ferential geometry, Vol.~2, Publish or Perish Press, 1979.

\bibitem{Steltsova09}
 Streltsova I.S.,
 $\mathbb{R}$-conformal invariants of curves,
 {\it Izv. Vyssh. Uchebn. Zaved. Mat.}  {\bf 53} (2009), no.~5, 67--69.

\bibitem{Stormark2000}
 Stormark O.,
 Lie's structural approach to PDE systems, {\it Encyclopedia of Mathematics and its Applications}, Vol.~80, Cambridge University Press, Cambridge, 2000.

\bibitem{Sulanke79}
 Sulanke R.,
 On \'E.~Cartan's method of moving frames, {\it Colloq. Math. Soc. Janos Bolyai}, Vol.~31, Dif\/ferential Geometry, Budapest, 1979.

\bibitem{Vassiliou1}
 Vassiliou P.,
 A constructive generalised Goursat normal form,
 {\it Differential Geom. Appl.}  {\bf 24} (2006), 332--350,
 \href{http://arxiv.org/abs/math.DG/0404377}{math.DG/0404377}.

\bibitem{Vassiliou2}
 Vassiliou P., Ef\/f\/icient construction of contact coordinates for partial prolongations,
 {\it Found. Comput. Math.} {\bf 6} (2006), 269--308,
\href{http://arxiv.org/abs/math.DG/0406234}{math.DG/0406234}.

\bibitem{Vessiot26}
 Vessiot E.,
 Sur l'int\'egration des faisceaux de transformations inf\/init\'esimales dans le cas o\`u, le degr\'e du faisceau \'etant $n$, celui du faisceau deriv\'ee est $n+1$,
 {\it Ann. Sci. \'Ecole Norm. Sup. (3)}  {\bf 45} (1928), 189--253.

\bibitem{Yamaguchi82}
 Yamaguchi K.,
 Contact geometry of higher order,
 {\it Japan. J. Math. (N.S.)} {\bf 8} (1982), 109--176.\\
Yamaguchi K.,
  Geometrization of jet bundles,
  {\it Hokkaido Math. J.} {\bf  12} (1983), 27--40.

\end{thebibliography}
\end{document}